\documentclass[12pt,a4paper,twoside,final,notitlepage, leqno]{article}
\usepackage[english]{babel}
\usepackage[T1]{fontenc}  
\usepackage{graphicx}
\usepackage{amsmath,amsthm,epsfig,amsfonts,bbm}  

\newcommand{\pequationdeb}{$$ \left\{ \begin{minipage}[c]{130mm}}
\newcommand{\pequationfin}{\end{minipage}
                           \right. $$}

\def \smb {{\scriptstyle \bullet }}
\newcommand{\monitem}{ \smallskip \noindent $\bullet$ \quad  } 
\newcommand{\moneq}{\vspace*{-6pt} \begin{equation} \displaystyle } 
\newcommand{\moneqstar}{\vspace*{-6pt} \begin{equation*} \displaystyle } 
\newcommand{\monendstar}{\vspace*{-6pt} \end{equation*}   }
\newcommand{\monend}{\vspace*{-6pt} \end{equation}   }
\newcommand{\beq}     {\begin{equation}}
\newcommand{\enq}     {\end{equation}}
\newcommand{\be}    {\begin{enumerate}}
\newcommand{\ee}    {\end{enumerate}}

\newcommand{\Bb}

\def\R{{\rm I}\! {\rm R}}

\def\section*#1{}

\def\resume{\if@twocolumn
\section*{R\'esum\'e}
\else \small
\quotation{\bf \it R\'esum\'e \rule[1mm]{1.5mm}{0.2mm}\vspace{0pt}}
\fi}
\def\endresume{\if@twocolumn\else\endquotation\fi}
\def\abstract{\if@twocolumn
\noindent\section*{{\bf Abstract}}
\else \small
\quotation{\noindent \bf {Abstract.} \rule[1mm]{1.5mm}{0.2mm}\vspace{0pt}}
\fi}
\def\endabstract{\if@twocolumn\else\endquotation\fi}

\usepackage{fancyhdr}
\renewcommand{\headrulewidth}{0pt}
\begin{document}

\fancypagestyle{plain}{ \fancyfoot{} \renewcommand{\footrulewidth}{0pt}}
\fancypagestyle{plain}{ \fancyhead{} \renewcommand{\headrulewidth}{0pt}} 

\title{\bf \LARGE    ~  \vspace{-.3 cm} ~\\   
 On triangular lattice Boltzmann schemes  ~\\ ~ \vspace{-.8cm}  ~\\ 
  for scalar problems     }

\author { { \large  Fran\c{c}ois Dubois~$^{ab}$ and   Pierre Lallemand~$^{c}$}     \\ ~\\ 
{\it  \small $^a$   Conservatoire National des Arts et M\'etiers, }  \\
{\it  \small  Department of Mathematics,  Paris, France. }  \\      
{\it \small  $^b$  Department of Mathematics, University  Paris-Sud,} \\
{\it \small B\^at. 425, F-91405 Orsay Cedex, France} \\ 
{\it \small francois.dubois@math.u-psud.fr} \\ 
{\it \small  $^c$   Beijing Computational Science Research Center,   Beijing Run Ze Jia Ye, China.} \\ %
{\it \small  pierre.lallemand1@free.fr}  \\  ~\\~\\ }  

\date{ 13 february   2012~\protect\footnote{~This contribution 
is published in {\it  Communications in Computational Physics}, 
 volume 13, number~3, pages~649-670, doi: 10.4208/cicp.381011.270112s, march 2013.
It is 
issued from a lecture  entitled ``D2T4 lattice Boltzmann scheme 
for scalar problems'' given on  monday 08 august  2011 at the conference  
``Discrete Simulations of  Fluid Dynamics'', 
 Fargo, North Dakota, USA.  }}

\maketitle

\begin{abstract} 
We propose to extend the d'Humi\`eres version of the lattice Boltzmann scheme to
triangular meshes. 
We use Bravais lattices or more general lattices with the property that 
the degree of each internal vertex is supposed to be constant. 
On such meshes, it is possible to define the lattice Boltzmann scheme as 
a discrete particle method, without need of finite volume formulation or 
Delaunay-Voronoi hypothesis for the lattice. We test this idea for the heat equation 
and perform an asymptotic analysis with the Taylor expansion method for two
schemes named D2T4 and D2T7. 
The results show a convergence up to second order accuracy and set new questions
concerning a possible super-convergence.
 $ $ \\[4mm]
   {\bf Keywords}: Laplacian operator, heat equation, d'Humi\`eres scheme, D2T4, D2T7.    
 $ $ \\[4mm]
   {\bf AMS classification}: 65-05, 65Q99, 82C20.   
\end{abstract}

\bigskip \bigskip  \newpage \noindent {\bf \large 1) \quad  Introduction}  

\fancyfoot[C]{\oldstylenums{\thepage}}
 

  \noindent 
The importance of extending the lattice Boltzmann scheme from 
square type regular meshes  to unstructured triangulations   has been 
recognized during  the last years of 20th   century \cite{Chen98, KSO99, PXDC99}.
In particular the ``volumetric formulation''  of 
Chen \cite{Chen98}  
makes a link with finite volumes, using control volumes  
around each vertex (the ``Inria cells'' \cite{Vi86}) 
  of  a finite element type triangulation.
This method is still under active development with the work of
Succi, Ubertini and co-workers 
\cite{PUS09, UBS03, UBS04}.
In a dual way,  van der Sman \cite{DE99,  DE2k,  vdS03, vdS04}  
uses rectangles and triangles as control volumes with a ``cell center'' type approach
in  Roache \cite{Ro72}   denomination. 
He has developed an approximation of diffusion equation with Delaunay-Voronoi
meshes for a BGK variant of the lattice Boltzmann scheme.

In a previous contribution \cite{DL08}, we have observed that for usual 
lattice Boltzmann schemes (as for example the well known D2Q9), 
several (two for D2Q9)  families of finite volumes 
are naturally associated with the scheme. 
As a consequence, we consider now the lattice Boltzmann scheme 
essentially as a ``particle'' method on a given ({\it a priori} fixed) mesh
with discrete velocities. 
Recall that the ``Particle In Cell'' 
 method has been first proposed  in 1964  by Harlow {\it et al.} 
\cite{HER64} and has been analyzed in the eighties by Beale and Majda \cite{BM82}, 
Raviart, Cottet and Mas Gallic  \cite{CG90,  GR87,  Ra85} among others. 
We remark  that this particle method does not suppose   {\it a priori} the 
 existence of a  given lattice. The surrounding cells
are recomputed at each time  step in order to make the particle interact.
Dynamic triangulation  is an alternative to the previous methodology. 
It has been developed recently by  Cianci, Klales, Love  and co-workers 
\cite{KCNML10, LC11} in the context of  lattice gas automata.

In this contribution, instead of adopting the  volumetric formulation
or a Delaunay-Voro\-noi hypothesis, 
we  develop the framework of  lattice Boltzmann schemes as a
variant  of the particle method. 
We propose an extension of   the approach of d'Humi\`eres \cite{DdH92}   
to  triangular meshes and we restrict this first tentative to scalar problems
like the heat equation without advection.

The outline of the contribution is the following. 
We first recall the classic D2T7 
 lattice  Boltzmann scheme in the next section. 
At this occasion, we put in evidence  a property of symmetry of 
Bravais lattices. 
It is possible to adapt the Taylor expansion analysis \cite{Du07} 
to this triangular lattice, with a diffusive scaling. This development is
presented in  Section~3 and applied to the D2T7 scheme. Several simulations
with the D2T7 lattice Boltzmann scheme for the heat equation are presented 
in  Section~4. In Section~5, we set the question of defining a 
discrete particle method on a finite element type triangular lattice. We propose a partial
answer when each vertex of the lattice has a constant number of neighbours. 
This framework is applied in Section~6 to define a D2T4 lattice Boltzmann scheme
for the heat equation. We repeat in Section~7 with this new scheme ``D2T4''
 the simulations presented in Section~4. This work validates
 the potential of applications of our proposal. The conclusion (section 8) serves also as 
a discussion concerning encountered difficulties.

 \newpage
\bigskip \bigskip   \noindent {\bf \large 2) \quad  D2T7 lattice  Boltzmann scheme}  


  \noindent   
We consider a Bravais lattice ${\cal L}$ connecting nodes labelled by the letter $x$ and
parametrized by a typical space scale  $\, \Delta x$. 
The  neighbour vertex number $j$ of the node $ \, x \in {\cal L} \, $ 
is denoted by $\, x_j \,$ and we set 
\moneq   \label{voisin-j}
x_j = x +  \xi_j \, \Delta x \, .  \monend  
For each $\, x \in {\cal L} \,$ and each direction $ \, \xi_j \,$ linking two vertices,
the ``opposite node'' with number  $\, \sigma(j) \,$ defined according to 
\moneq   \label{voisin-sigma-j}
x_{\sigma(j)} \equiv  x -  \xi_j \, \Delta x \,,\qquad 
\xi_j + \xi_{\sigma(j)} \equiv  0 \,   \monend  
is  also  a vertex of the lattice  ${\cal L}$ 
({\it i.e.} $ \, x_{\sigma(j)} \in {\cal   L}$).  
In the following, we emphasize this property satisfied by Bravais lattices
and qualify it as a symmetric property.  
%
Most ``DdQq'' schemes  (with a notation introduced by Qian {\it et al.} \cite{QHL92})
 presented in the literature use a Bravais lattice. 
%
This symmetry property is also mandatory {\it e.g.} to define ``two relaxation times'' 
lattice Boltzmann schemes as proposed by Ginzburg {\it et al.}  \cite{GVH08}.

The lattice Boltzmann scheme with multiple relaxation times
is defined in a classical manner. 
Consider  a vertex   $x$  that belongs to the lattice $ {\cal L}$.
Then  the    $ j ^{\rm o} $ direction of propagation is defined
with a vector  $ \xi_j $   and   $ \, \xi_j  \in  {\cal V} $,
set of directions that define the  vicinity of the vertex $x$. 
The    $ \, j ^{\rm o} \, $   density of particles
 at vertex $x$ and time $t$ is denoted by $  \, f_j(x, \, t) $.
After a local step of relaxation, the 
  $ \, j ^{\rm o} \, $   density of particles 
is named  $  \, f_j^*(x, \, t) $. 
Because a Bravais lattice $\,  {\cal L} \,$ is symmetric, 
 the  neighbouring vertex $\, x_{\sigma(j)} \,$ defined in 
(\ref{voisin-sigma-j}) in the direction {\bf opposite} 
 to the  $ j ^{\rm o} $ direction of  propagation 
  belongs to the lattice $  {\cal L}$. 
The  lattice Boltzmann scheme can be completely defined:  
\moneq   \label{LB-scheme} 
 f_j(x, \, t+\Delta t) =   f_j^*(x - \xi_j \, \Delta x , \, t) \, .
\monend 
Moreover the basic iteration (\ref{LB-scheme}) 
of a lattice Boltzmann scheme supposes explicitly that the lattice is symmetric, as
 illustrated  in Figure 1 (left). 

\smallskip    
\centerline { \includegraphics [width=.11 \textwidth] {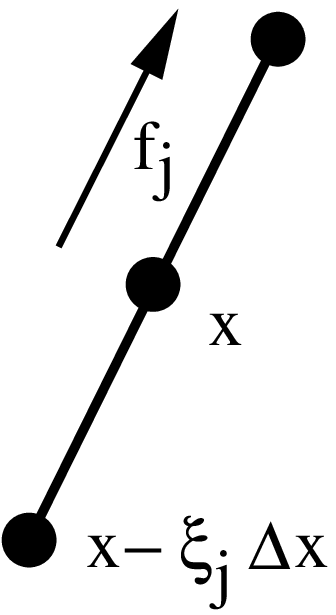}
\qquad \qquad \qquad  \qquad \qquad 
\includegraphics[width=.19 \textwidth] {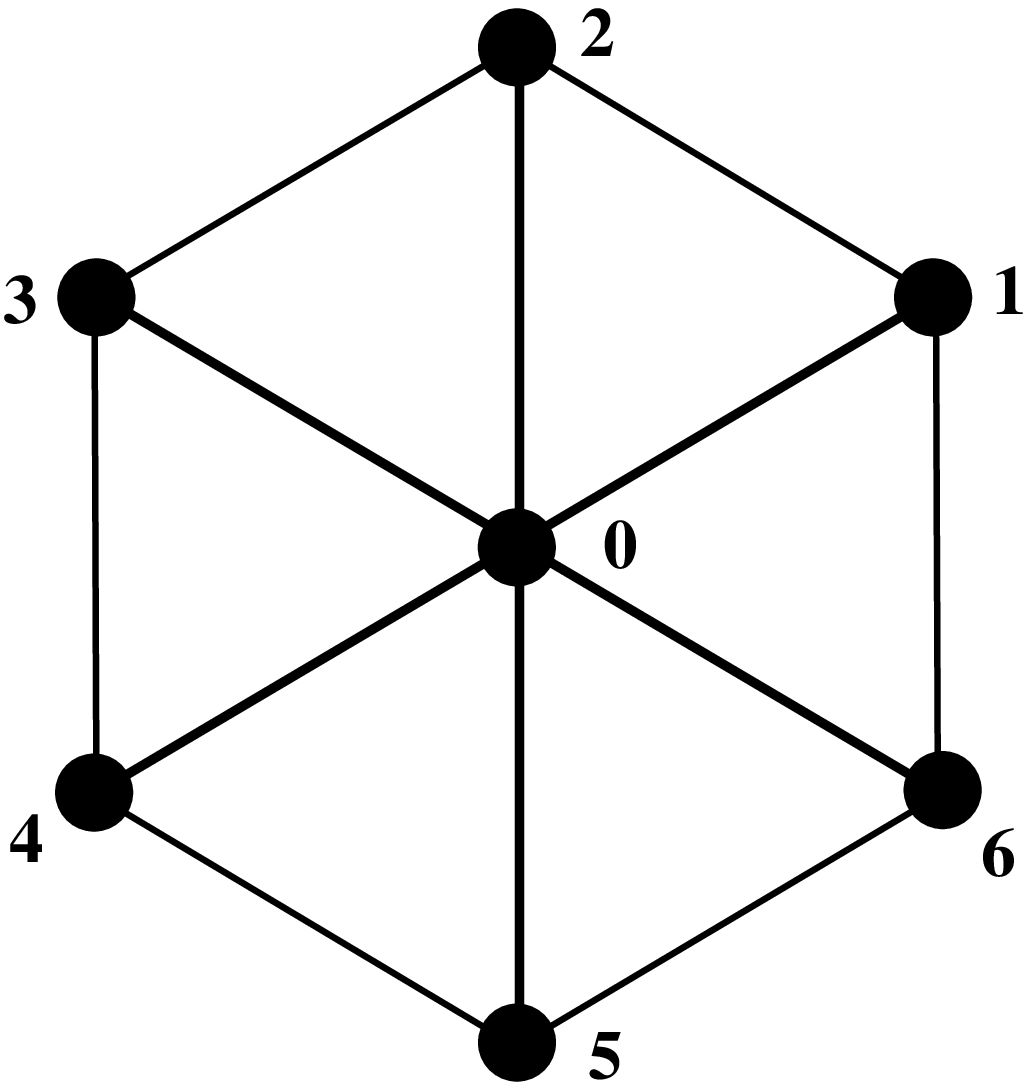} }  

\noindent  {\bf Figure 1}. \quad 
Typical stencil of a lattice Boltzmann scheme for a Bravais lattice (left)~;
 both opposite directions $\, \xi_j \,$  and  $\, -\xi_j \,$  connect 
two vertices of the   mesh.   Local numbering of the  six neighbours  (right)
of the  D2T7 lattice Boltzmann scheme on triangles.   
\smallskip \smallskip 


\smallskip    \noindent 
The D2T7 lattice Boltzmann uses equilateral triangles as suggested
 by Frisch,  Hasslacher and Pomeau in 1986    \cite{fhp} 
 in the context of   lattice gas automata. 
We precise   the parameters that we have to consider. 
A vertex $x$  has a total of six  neighbours (seven including itself) 
$ \, \xi_j $ $(j=0, \dots , \, 6$) as in Figure~1 (right).  Following  
d'Humi\`eres approach  \cite{DdH92}  we introduce 
moments    $ \, m_k \,$  as linear functions of the particle distribution
$ \, f \,$:  
\moneq   \label{f-to-m} 
 m_k = \sum_j M_{kj} \, f_j   \, .
\monend 
We restrict our study to the  
simple case of only one conservation (thermal  problem).
Following \cite{LL00}, we introduce a 
family $\, {\cal P} \,$ of polynomials  $\, p_k \, $ for $ \, k=0, \cdots , \, 6$:  
\moneq   \label{polynomes} 
 {\cal P} = \big\{  1, \, X ,\,  Y ,\, X^2 + Y^2  ,\,   
{{4}\over{\sqrt{3}}}\, X \, Y ,\,  2 \, ( X^2 - Y^2) ,\,  3\, Y - 4 \, Y^3 \big\} \, . 
\monend 
The coefficients of the matrix $\, M \,$ introduced at relation 
(\ref{f-to-m}) are simply given by a nodal value in the velocity space: 
\moneq   \label{matrice-M} 
M_{kj} = p_k(\xi_j) \,, \qquad 0 \leq \, j , \, k \, \leq 6 \, . 
\monend 
We remark that    $ \, M_{\alpha j} = \xi_j^{\alpha}  \, $ 
for $\, \alpha = 1, \, 2. \, $ 
We have only one conserved moment 
$ \,  \rho \equiv m_0 = m_0^{\rm eq}= m_0^*  = \sum_j f_j \, $  
and the other moments at equilibrium follow the relations 
$ \,\,  m_1^{\rm eq} = m_2^{\rm eq} = 0 , \,$
$ \,  m_3^{\rm eq} = a_3 \, \rho  ,\, $
$ \,\,  m_4^{\rm eq} = m_5^{\rm eq} = m_6^{\rm eq} = 0 .\, $
The relaxation of moments out of equilibrium is also very simple: 
\moneq   \label{m-star}  
 m_k^* = m_k + s_k \, (m_k^{\rm eq} - m_k) \,, \quad  k = 1, \cdots , \, 6 \, . 
\monend 
with $ \, s_1 = s_2 \, $ and $ \, s_4 = s_5 \, $ to enforce isotropy.

\bigskip \bigskip   \noindent {\bf \large 3) \quad  Taylor expansion with diffusive scaling}  


 \noindent    
We can analyse the D2T7 lattice Boltzmann scheme with the 
 Taylor expansion method    \cite{Du07}.
We consider one time step of iteration (\ref{LB-scheme}) 
and we replace the particle distribution in the right hand side by the moments after
 relaxation: 

\smallskip \noindent 
$ \displaystyle     f_j(x, \, t+\Delta t) =  \sum_\ell M^{-1}_{j \ell} \,
 m_\ell^*(x - \xi_j \, \Delta x , \, t) . \,\, $

\smallskip \noindent 
In consequence, we have the formal expansion in the moment space

\smallskip \noindent 
$ \displaystyle   m_k(x, \, t+\Delta t) = \sum_{ j \ell} M_{k j}  \,  M^{-1}_{j \ell} 
\, \, m_\ell^*(x - \xi_j \, \Delta x , \, t) $

\smallskip \noindent 
 $ \displaystyle  \qquad \qquad \qquad  \,\,  \,  = \,   \sum_{ j \ell} M_{k j}  \,  M^{-1}_{j \ell} \, \big[ 
 m_\ell^*(x    , \, t) -  \xi_j^\alpha \, \Delta x \,  \partial_\alpha m_\ell^*
+ {\rm O}(\Delta x^2) \big]   $  

\smallskip \noindent 
  $ \displaystyle  \qquad \qquad \qquad    \,\,  \,  = \,   m_k^*  -  \Delta x \,  \sum_{\ell} 
\big( \sum_{j} M_{k j} \, \, \xi_j^\alpha \,\,    M^{-1}_{j \ell} \big) \, 
 \,  \partial_\alpha m_\ell^* + {\rm O}(\Delta x^2)   .\,  $  

\smallskip \noindent 
We introduced  the    momentum-velocity tensor introduced in \cite{Du07}: 
$   \displaystyle 
 \Lambda_{k p}^\ell \equiv   \sum_{j} M_{k j} \, \,   M_{p  j} \, \,    M^{-1}_{j \ell}
 \, . $
Then we have up to third order accuracy 

\smallskip \noindent 
$  \displaystyle   m_k(x, \, t+\Delta t) = 
 m_k^*  -  \Delta x \,  \Lambda_{k \alpha}^\ell \, 
  \partial_\alpha m_\ell^* + {1\over2} \,  \Delta x^2  \,\,   
\Lambda_{k \alpha}^p \,  \Lambda_{p \beta}^\ell \,    
\partial_\alpha  \partial_\beta m_\ell^* + {\rm O}(\Delta x^3) \,  $

\smallskip \noindent 
and using the relaxation step (\ref{m-star}), 

\smallskip \noindent 
 $  \displaystyle m^*_k(x, \, t+\Delta t) = 
 m_k^{\rm eq}   -  {\Delta x}\, {{1-s_k}\over{s_k}} \,  \Lambda_{k \alpha}^\ell \,  
 \partial_\alpha m_\ell^{\rm eq}  + {\rm O}(\Delta x^2) . \, $

\smallskip \noindent 
We adopt the so-called ``diffusive scaling'' proposed  
initially for rarefied flows by Sone \cite {So69} 
(see  an explicit derivation for lattice Boltzmann schemes  {\it e.g.}
in Junk {\it et al.}  \cite{JKL05})  
\moneq   \label{diffusive-scaling} 
\Delta t \equiv {{\Delta x^2}\over{\zeta}}   \,     \monend
where $\, \zeta \,$ is a constant for homogeneity of dimensions.  
We  add some advection term by enforcing   the relations according to 
$ \,\,  m_1^{\rm eq} =  u \,  {{\Delta x}\over{\zeta}} \,\, $ and 
$ \,\,  m_2^{\rm eq} =  v \,  {{\Delta x}\over{\zeta}}  . \,$
After some pages of formal calculus, following the method presented in details
in  \cite{DL09},  we obtain 
the equivalent partial differential equation :
\moneq   \label{equiv-pde-order4}  
  {{\partial \rho}\over{\partial t}} + u \,  {{\partial \rho}\over{\partial x}} 
+ v \,  {{\partial \rho}\over{\partial x}} - \mu \,  \Delta \rho  \, = \,  
 \Theta \, \Delta x^2   \,  \,    \Delta^2  \rho   \, + \, 
 \Delta x^4 \, A_6 \, \rho +   {\rm O}(\Delta x^6)   \,.  \monend 
Up to second order accuracy, we have an approximation of the heat equation  
with  a  diffusivity coefficient $ \, \mu \,$   given according to 
\moneqstar   \mu =  {1 \over 2 } \, \zeta \, a_3 \, \sigma_1 .    \monendstar
The  coefficients $ \, \sigma_k \, $  for the nonconserved moments  are
given by the  H\'enon's relation \cite{He87} 
$\,\,  \sigma_k \equiv  {{1}\over{s_k}} - {1 \over 2 }  . \, $ 
The coefficient $ \, \Theta \,$ in front of the fourth order term 
in (\ref{equiv-pde-order4}) is explicited as follows for $\, u=v=0 \,$:
\moneq   \label{Theta} 
\Theta = - {{1}\over{16}} \, \sigma_1 \, a_3 \, \zeta \, \Big( (1 - a_3) \, 
\big( 1 - 4 \, \sigma_1 \, \sigma_3 \big)  - 2 \,  \sigma_1 \, \sigma_4 + 4 \, a_3 \,
\sigma_1^2 \Big)  \, .
\monend  
In the relation (\ref{equiv-pde-order4}), $\, A_6 \,$ is a sixth order operator. 
The development of the other moments can also be achieved.  
In particular, we have 
$ \,\,  m_\alpha = m_\alpha^{\rm eq}  - {{a_3}\over{2 \, s_1}} \, 
\Delta x \, \partial_\alpha  \rho  +   {\rm O}(\Delta x^2)   . \, $ 
%

\smallskip \monitem  
``Second order'',   ``quartic'' and ``hexahedric'' coefficients 

 \noindent
We have chosen    the following numerical values  
$ \,  \zeta = 1    \,, \,   a_3  = {1\over4}  \,, \,   s_1  = 0.8  \,$  
compatible  with a diffusivity coefficient 
  $   \, \mu =     0.09375 $. In these  conditions, 
the D2T7 lattice Boltzmann scheme is formally equivalent to the heat equation 
 up to order 2 ({\it id est}, due to (\ref{equiv-pde-order4}) and (\ref{Theta}), 
  $ \, \Theta \neq 0 $ and $ A_6 \neq 0 $) 
when using to fix the ideas the following ``second order''  coefficients (given here
with 15 decimals for a possible implementation): 
\moneq   \label{coefs-order2} 
 s_3    =    1.428571428571428   ,\,\, 
 s_4 = s_5 = 0.481927710843373    ,\,\, 
 s_6       = 0.476190476190476   .   \monend  
%
%
With the following choice of ``quartic'' relaxation coefficients  
\moneq   \label{coefs-order4} 
 s_3    =    1.428571428571428 ,\,\, 
 s_4 = s_5 = 0.930232558139534  ,\,\, 
 s_6       = 0.526315789473684   ,    \monend  
%
%
the D2T7 lattice Boltzmann scheme is formally of the  order 4 ({\it id est} 
$ \, \Theta =  0 $ and $ A_6 \neq 0 $). Last but not least, 
we can impose  $ \, \Theta \equiv 0 $ and $ A_6 \equiv 0 $ 
and the D2T7  scheme is  of  order 6. The 
``hexahedric''  coefficients  can be taken as follows: 
\moneq   \label{coefs-order6} 
 s_3    =    1.086117521785847  ,\,\, 
 s_4 = s_5 = 1.344205296559553  , \,\, 
 s_6       = 0.647305233773416  .    \monend  
%
%

\bigskip \bigskip   \noindent {\bf \large 4) \quad Diffusion simulations   with the    D2T7 scheme}  


\noindent  
We have done several simulations: a ``one point'' periodic analysis, a
numerical evaluation of the  modes for a periodic pipe and a rectangle, the computation of
harmonic functions by time asymptotics of the heat equation, the dissipation of a 
triangular Dirichlet  mode and the direct numerical computation of triangular
 Dirichlet  modes. 

\newpage 
\smallskip \monitem     One point periodic analysis     

\noindent 
The one point  analysis   can be conducted as follows. We start from the 
iteration (\ref{LB-scheme}) of a lattice Boltzmann scheme. We suppose that 
the particle field for the neighbouring points 
of vertex $x$ satisfy the following  periodicity condition :  
\moneq   \label{periodic} 
f_j \big( x - \xi_j \, \Delta x , \, t) \,=\,  \exp 
\big( -i \,   {\rm \bf  k} \, \smb \,  \xi_j  \, \Delta x \big) \, \,  f_j \big( x , \, t)
\monend  
for some wave vector $ \, {\rm \bf  k} = (k \cos \theta ,\, k \, \sin \theta ) .$ 
From (\ref{periodic}), the evaluation of the right hand side of  (\ref{LB-scheme}) 
is easy in the context of the d'Humi\`eres  version of the lattice Boltzmann scheme. 
The state vector $f$ is then solution of an eigenvalue problem of small dimension $q$
for a general lattice Boltzmann problem with $q$ velocities. In the D2T7 case
for thermal problems, we obtain six eigenvalues $ \, \lambda_\ell \simeq 1 - s_\ell \,$
for $ \, \ell \geq 1 \,$ and one physical eigenvalue $\, \lambda( {\rm \bf  k} )
\simeq 1 -\mu \, k^2 $.  This eigenvalue has a  real meaning for  applications
to macroscopic physics. 
A numerical diffusivity $ \, \mu_{\rm num} \equiv (1-\lambda( {\rm \bf  k} )) / k^2 \, $
can be extracted from the previous relation. In Figure 2, we have plotted the error
$ \, \epsilon \equiv \, \mid \mu -  \mu_{\rm num} \mid \,$ as a function of the modulus of
the wave vector. 
With the three versions of the D2T7 scheme detailed in 
(\ref{coefs-order2}), (\ref{coefs-order4}) and  (\ref{coefs-order6}),
the errors for the diffusivity have an order of convergence directly predicted by the Taylor
expansion analysis.

\smallskip    
\centerline { \includegraphics[width=.47 \textwidth] {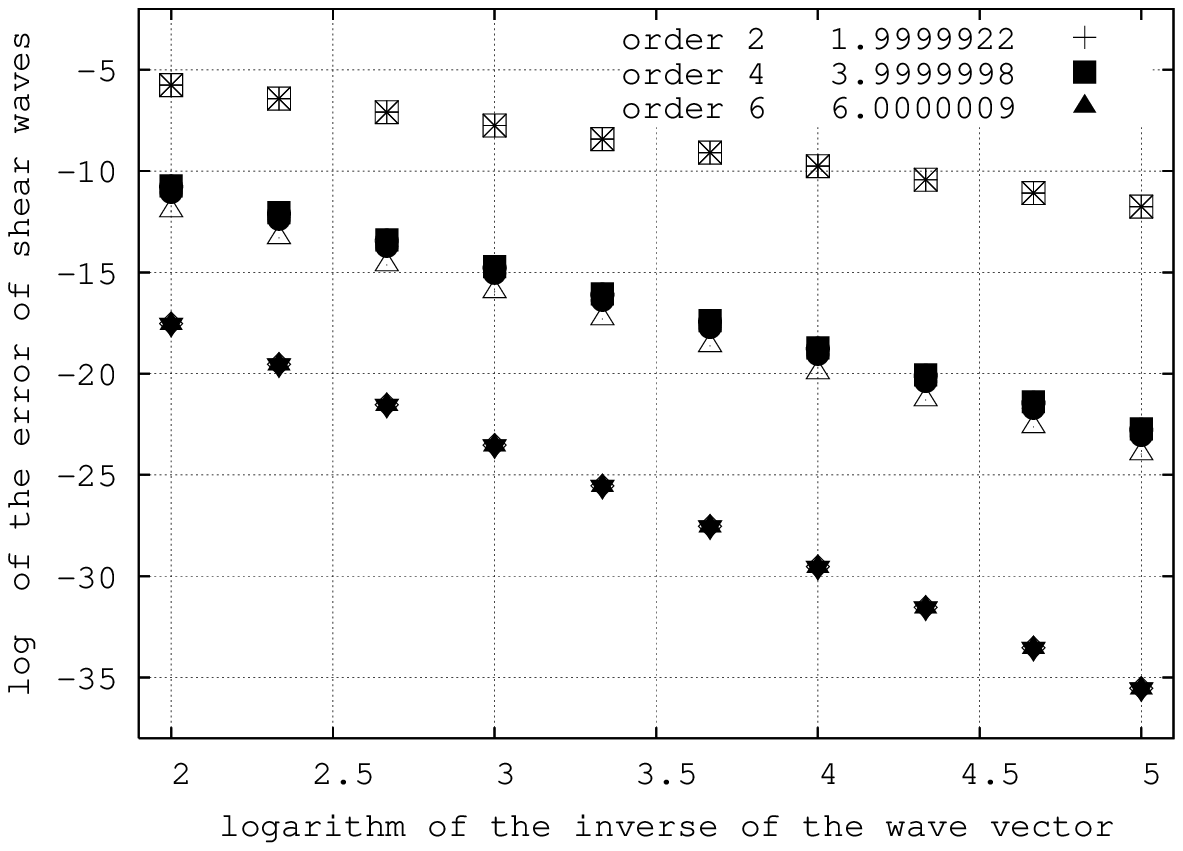}  } 

 \noindent  {\bf Figure 2}. \quad 
D2T7 lattice Boltzmann scheme for the heat equation. One point periodic analysis. 
Error $ \epsilon \, \equiv \,   \mid \mu - \mu_{\rm num} \mid \, $  between  
exact and numerical  diffusivities.
Note that the exact dispersion equation can be obtained, and when solved by successive
approximations in powers of $k$, it leads to the same results, obtained from the 
successive equivalent equations.    
\smallskip \smallskip 

 \smallskip   \monitem       Periodic pipe and rectangle 

\noindent 
The analysis for a periodic pipe is conducted by following the same ideas. 
A D2T7 lattice Boltzmann solver is considered  on a simple geometry of 
  $ \, nx \equiv 96  \,$ by $\, ny \equiv 4 \, $ mesh points. The unknown is now a vector 
$ \, f \in \R^{\, 7 \,  nx  \, ny} . \, $  The iteration of the scheme defines a linear
operator $\, A \,$ and the first eigenvalue  of this operator is determined thanks
to an Arnoldi algorithm \cite{Ar51}. 
The first eigenvalue $ \, \lambda \equiv 1 \,$ corresponds to the conservation 
of mass in the whole domain, including boundary conditions. 
The second eigenvalue $\, \lambda_0 \,$ 
corresponds to the smallest wave vector compatible with the computational domain. It 
is compared with the modulus of the wave vector to evaluate a numerical diffusivity 
$ \, \mu_{\rm num} =  (1-\lambda_0) / k^2 \,$ as previously. 
The different errors $ \, \epsilon \equiv \, \mid \mu -  \mu_{\rm num} \mid \,$
are presented in Figure~3 (left). 
The first two versions (\ref{coefs-order2})(\ref{coefs-order4}) of orders
two and four  present a coherent numerical convergence. 
The results are not so clear with the sixth order tuning of the parameters.
It seems to be due to the round-off errors for this study involving
three orders of magnitude for wave vector. 
 
The analysis is analogous for a rectangle 
  $ \, nx \equiv 36 \,$ by 
$ \, ny \equiv 52 \, $ mesh points.    
The results are depicted in Figure~3 (right). 
The lattice Boltzmann scheme has a coherent order of convergence for the 
``second order'' and ``fourth order'' versions of the scheme. 
The ``sixth order'' scheme exhibits now an error
numerically evaluated as only fifth order accurate. 
This fact seems again to due to round-off errors in the Arnoldi process \cite{Ar51}.

\smallskip      
\centerline {   \includegraphics  [width=.47 \textwidth]   {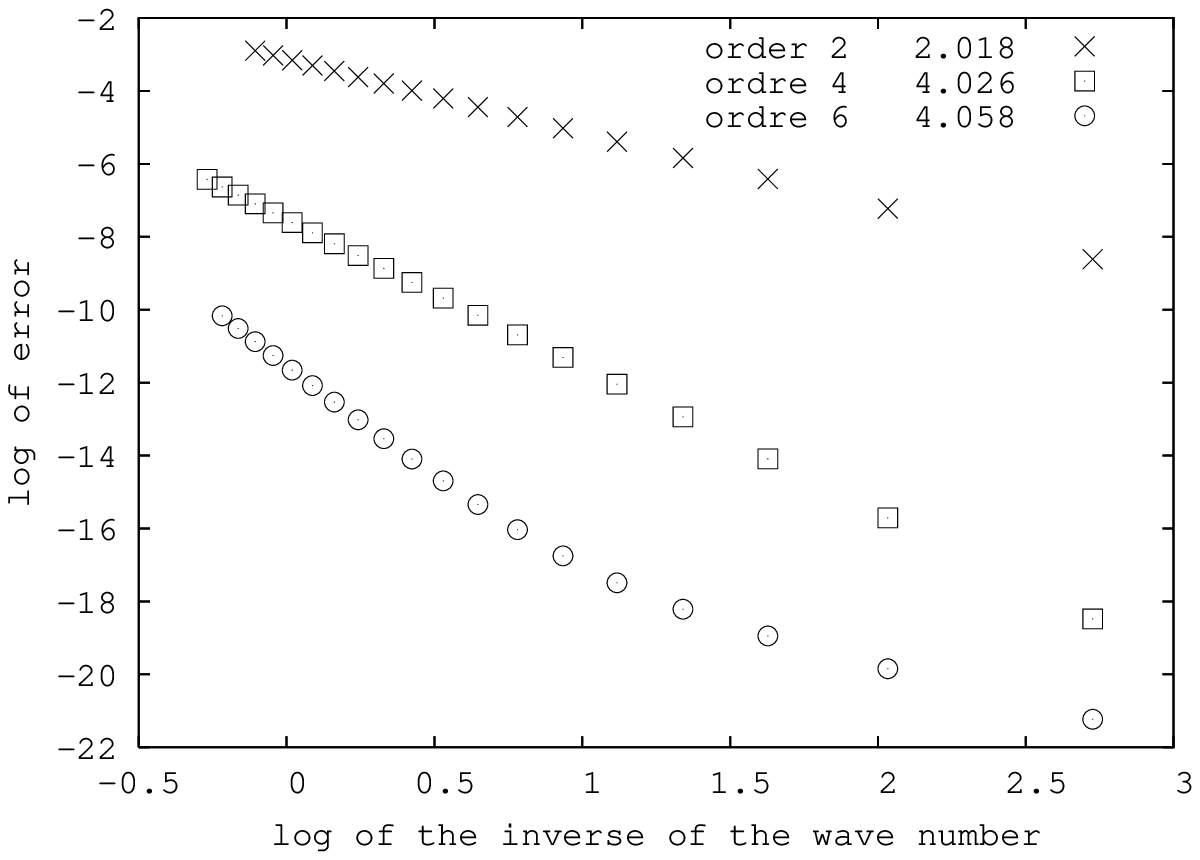} 
\quad 
\includegraphics  [width=.47 \textwidth] {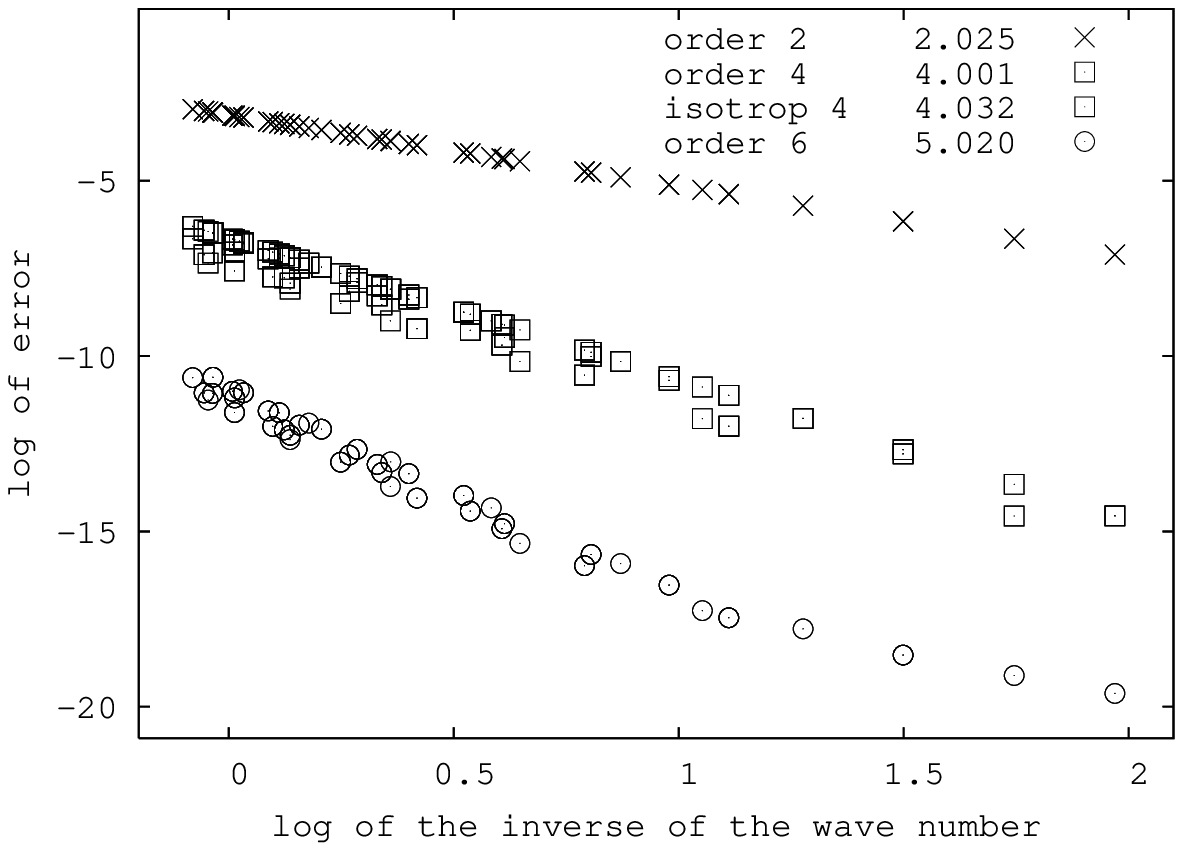} } 
 
\noindent  {\bf Figure 3}. \quad  
D2T7  lattice Boltzmann scheme for the heat equation. 
Periodic modes for a pipe with $ nx = 96 $  and $ ny = 4 $ mesh points (left)
an   periodic modes for a rectangle
of $ nx = 36 $  by  $ ny = 52 $     points (right). 
Error $ \epsilon \, \equiv \,   \mid \mu -  \mu_{\rm num}  \mid \, $  
between  exact and numerical diffusivities. The hexahedric predicted coefficients  
 define a fourth or fifth order scheme.    
\smallskip \smallskip 

\smallskip       
\centerline {  \includegraphics [width=.35 \textwidth, height=.42 \textwidth] 
{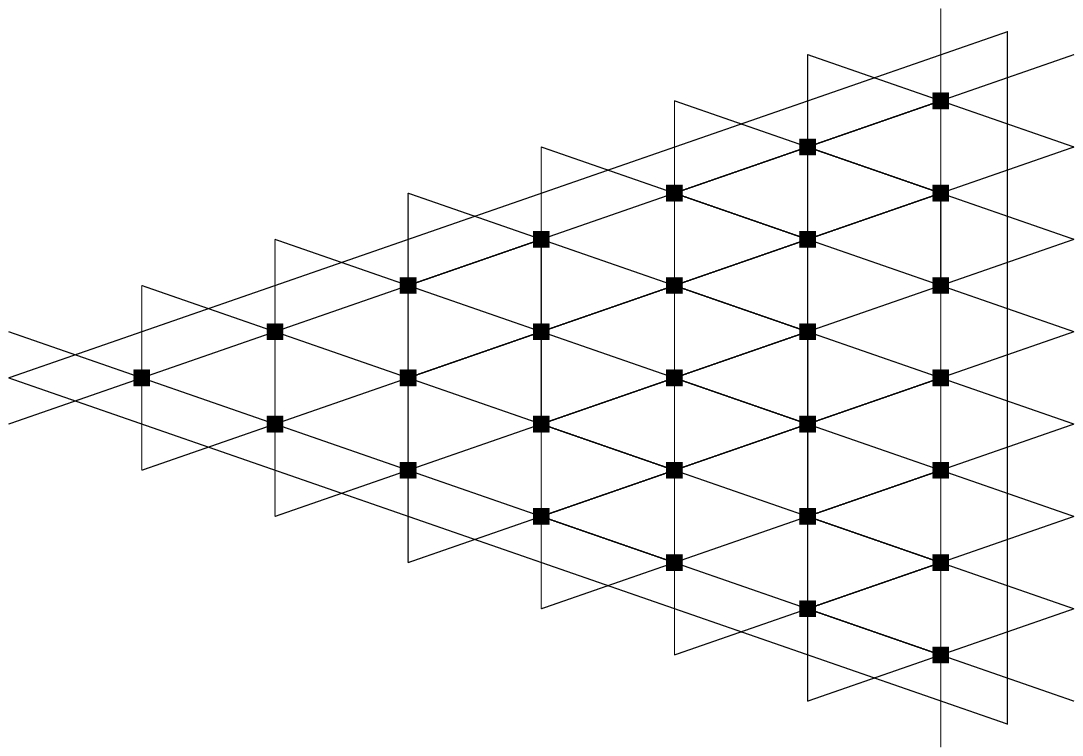} \qquad  \qquad 
 \includegraphics[width=.37 \textwidth] {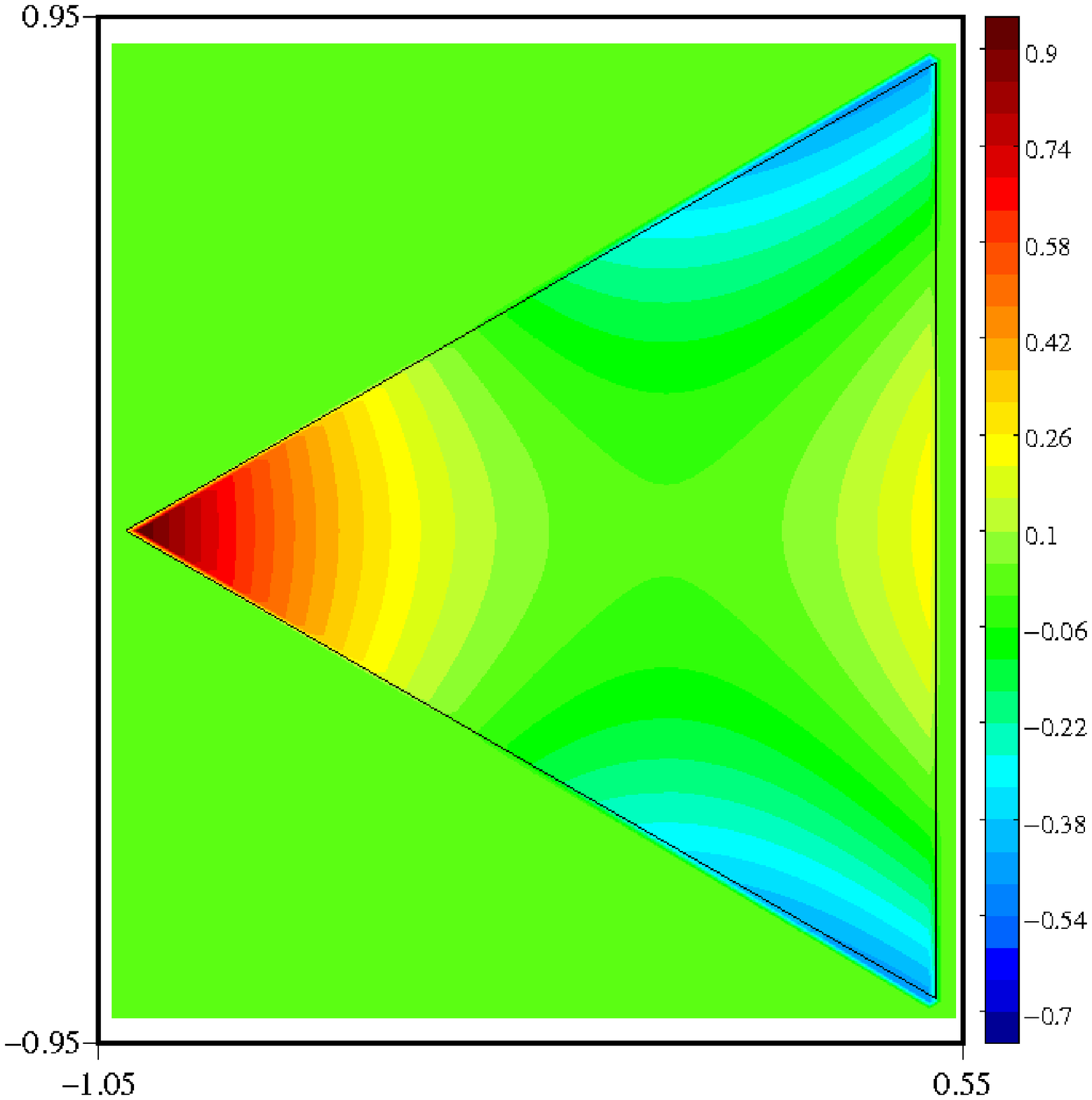} }

\noindent  {\bf Figure 4}. \quad 
Typical two-dimensional mesh for a D2T7 computation on a triangle (left).  
Two-dimensional computation of the harmonic function  $ \, p_H(x,\, y)  = x^2 - y^2 \, $  
on a triangle with the D2T7 lattice Boltzmann scheme (right). The iso-contours are 
composed by discrete hyperbolas.   
\smallskip \smallskip 

\newpage 
\smallskip  \monitem 
   Harmonic polynomials on a triangle     

\noindent
We have developed a D2T7 solver for a triangular geometry (see a typical mesh in
Figure~4) The initial condition is {\it a priori} identically null. We determine the
numerical boundary conditions compatible with a polynomial expression 
$ \,  p_H(x,\, y)  \equiv  x^2 - y^2 \,  $  on the boundary
with an ``anti-bounce-back'' version of the algorithm of Bouzidi {\it et al.} \cite{bfl01}. 
The computation converges in time towards the harmonic function introduced
above. 
We present in Figure~4 (right) the numerical result 
$ \,  \rho(x, \, y) \simeq  p_H(x,\, y) \, $ when we use $ \, 61 \,$ 
points on the edge of the triangle (that corresponds to a total of 
1891  vertices for the entire mesh).

\smallskip     
\centerline {   \includegraphics [width=.32 \textwidth] {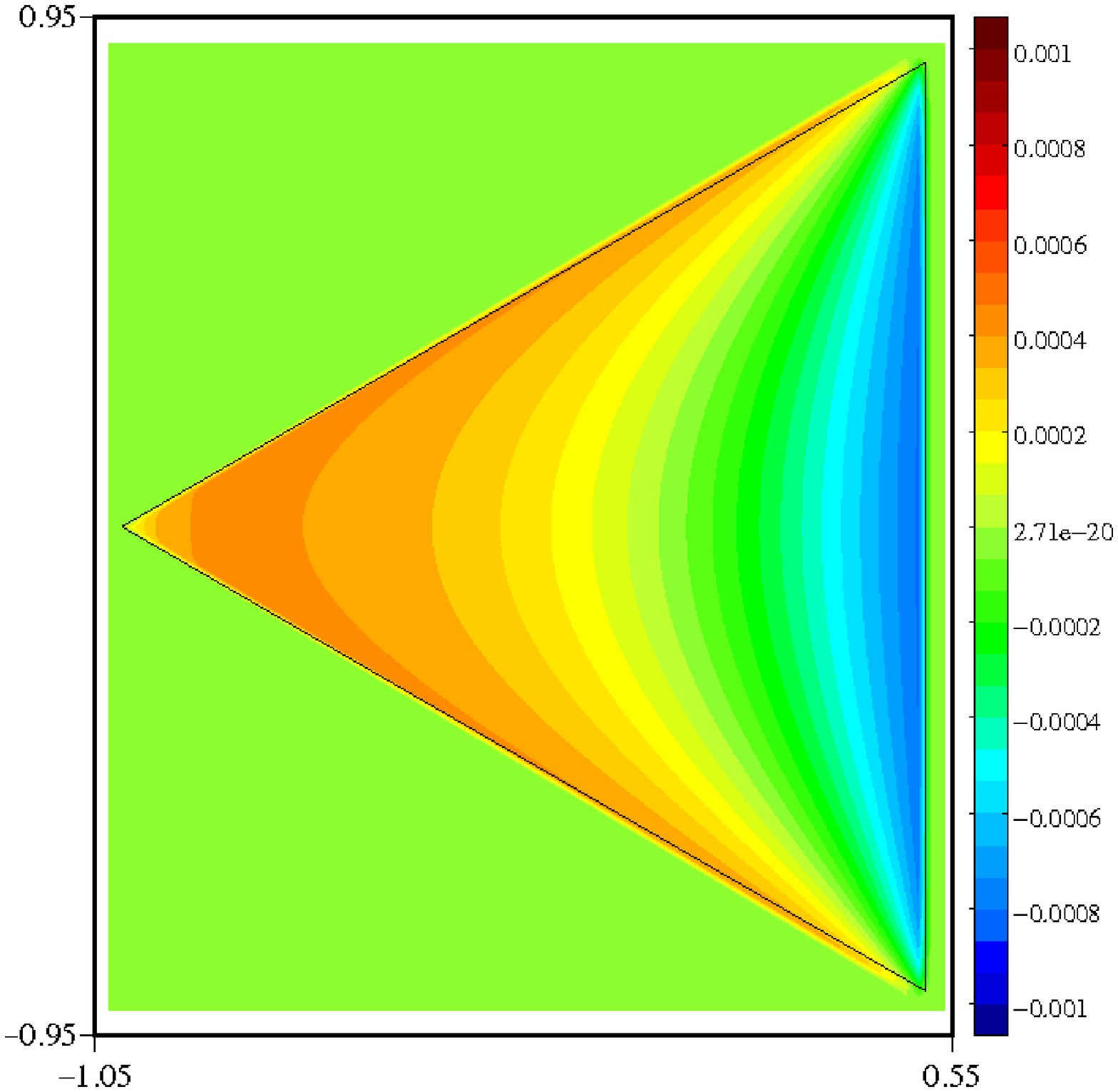} 
\includegraphics [width=.32 \textwidth] {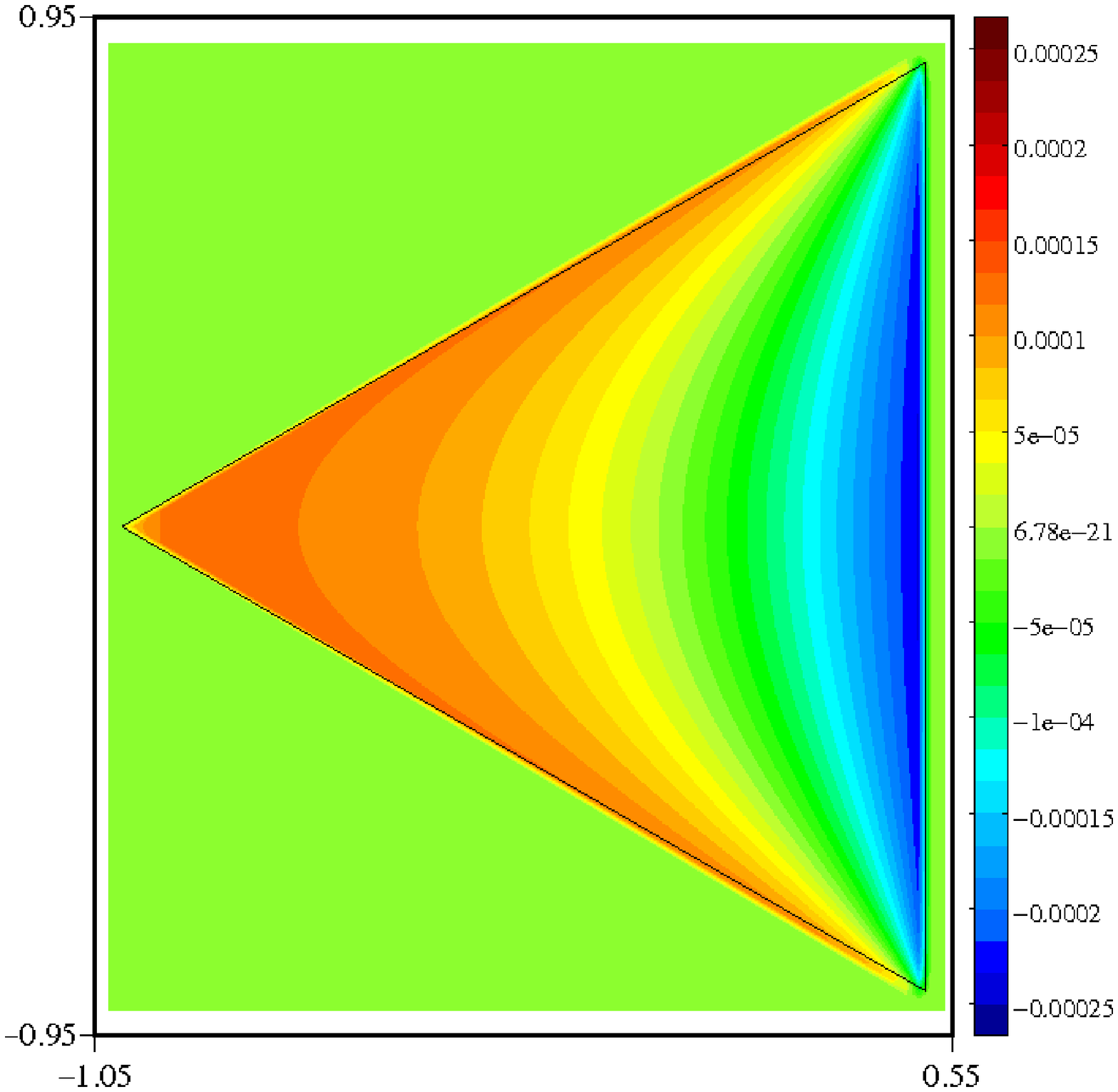} 
\includegraphics [width=.32 \textwidth] {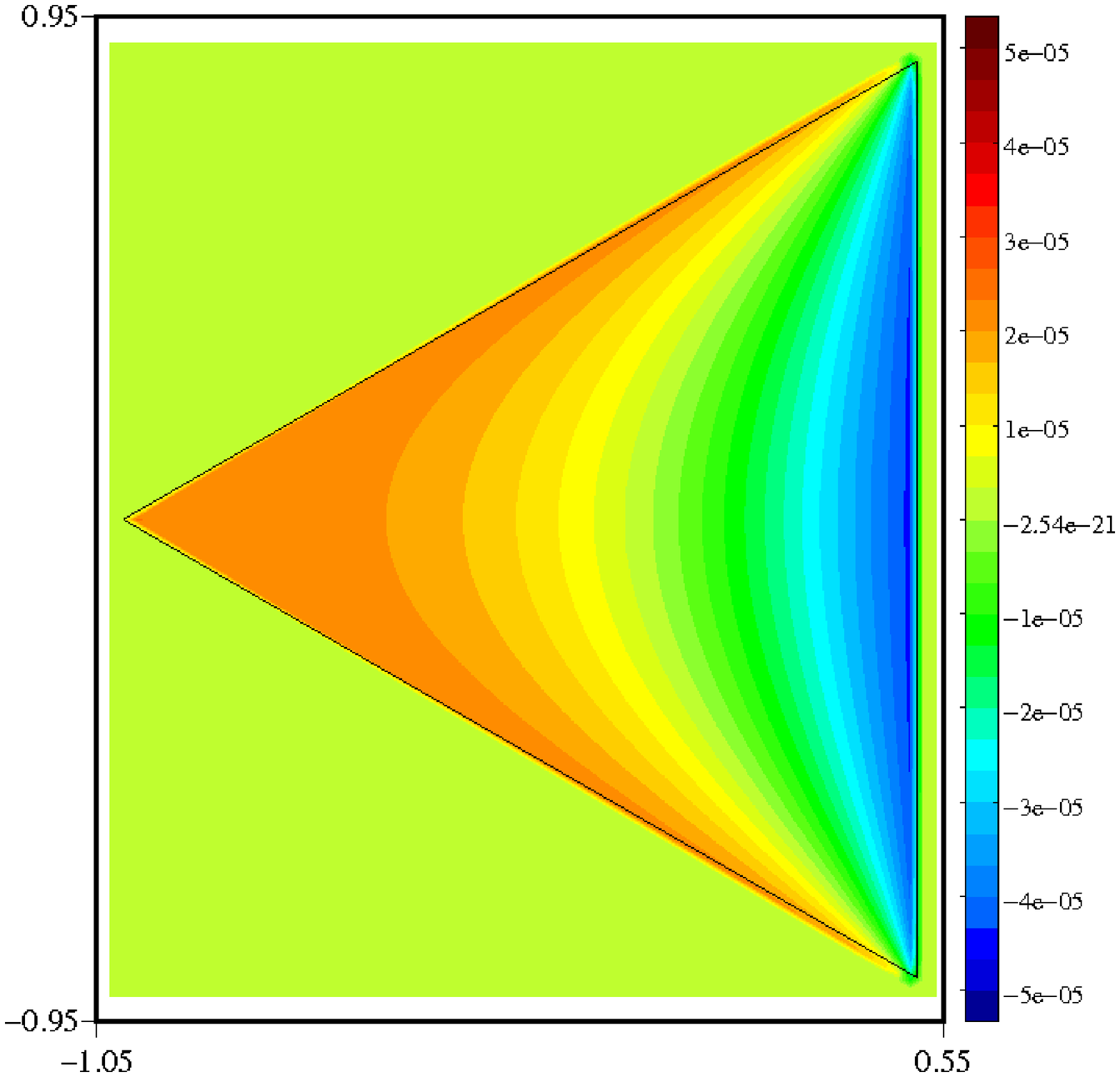} }

\noindent  {\bf Figure 5}.  \quad 
Two-dimensional computation of the harmonic function  $ \, p_H(x,\, y)  = x^2 - y^2 \, $  
on a triangle with the D2T7 lattice Boltzmann scheme. Iso-contours of the  errors 
for three sets  of parameters presented at relations (\ref{coefs-order2}), 
 (\ref{coefs-order4}) and  (\ref{coefs-order6}). 
 Negative values are in blue and positive ones in
red. The maximal errors are equal to 
 $ 8.14 \,\,  10^{-4} $ (left), $ 2.36 \,\,  10^{-4} $ (middle) and 
 $ 4.47 \,\,  10^{-5} $ (right) when using a D2T7 scheme with formal order of 
2  (left), 4  (middle) and 6 (right).  
\smallskip \smallskip 
 
\smallskip         
\centerline { \includegraphics[width=.47 \textwidth] {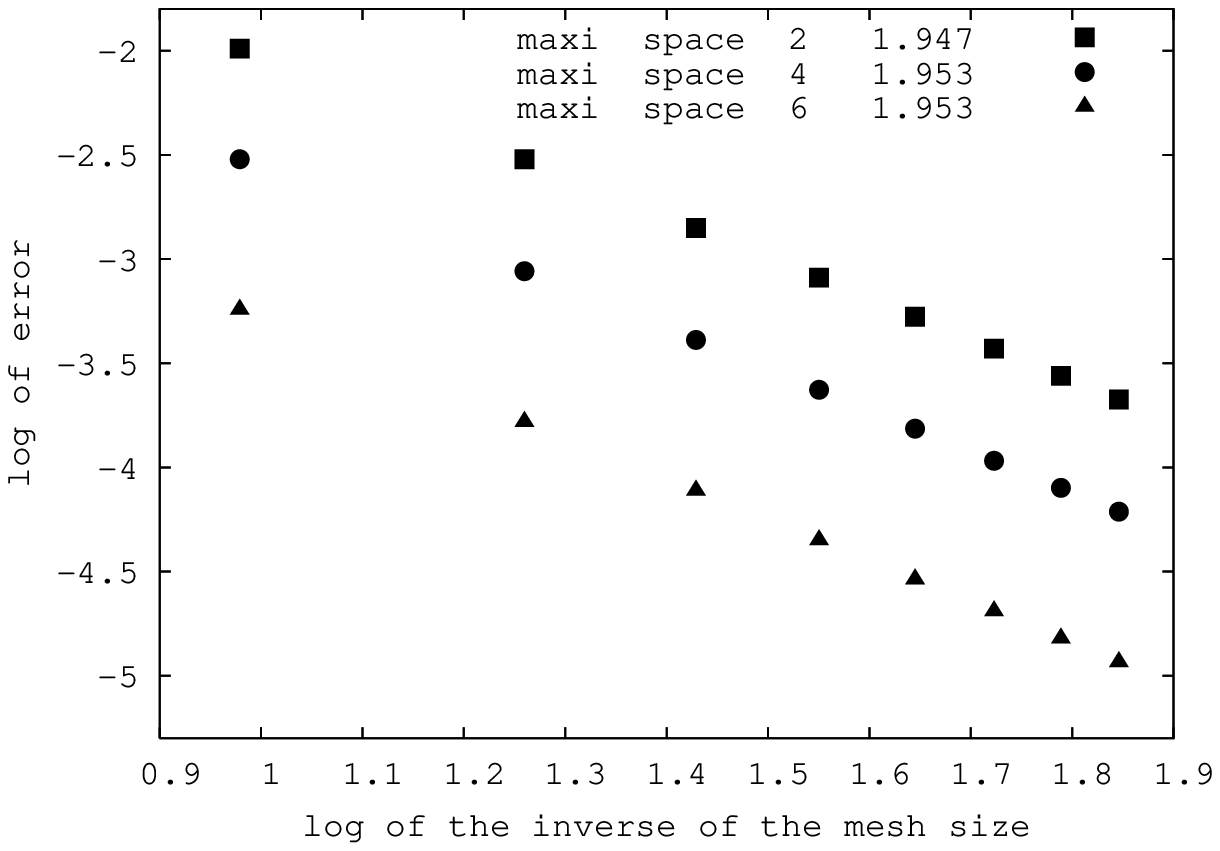}  \quad 
 \includegraphics[width=.47 \textwidth] {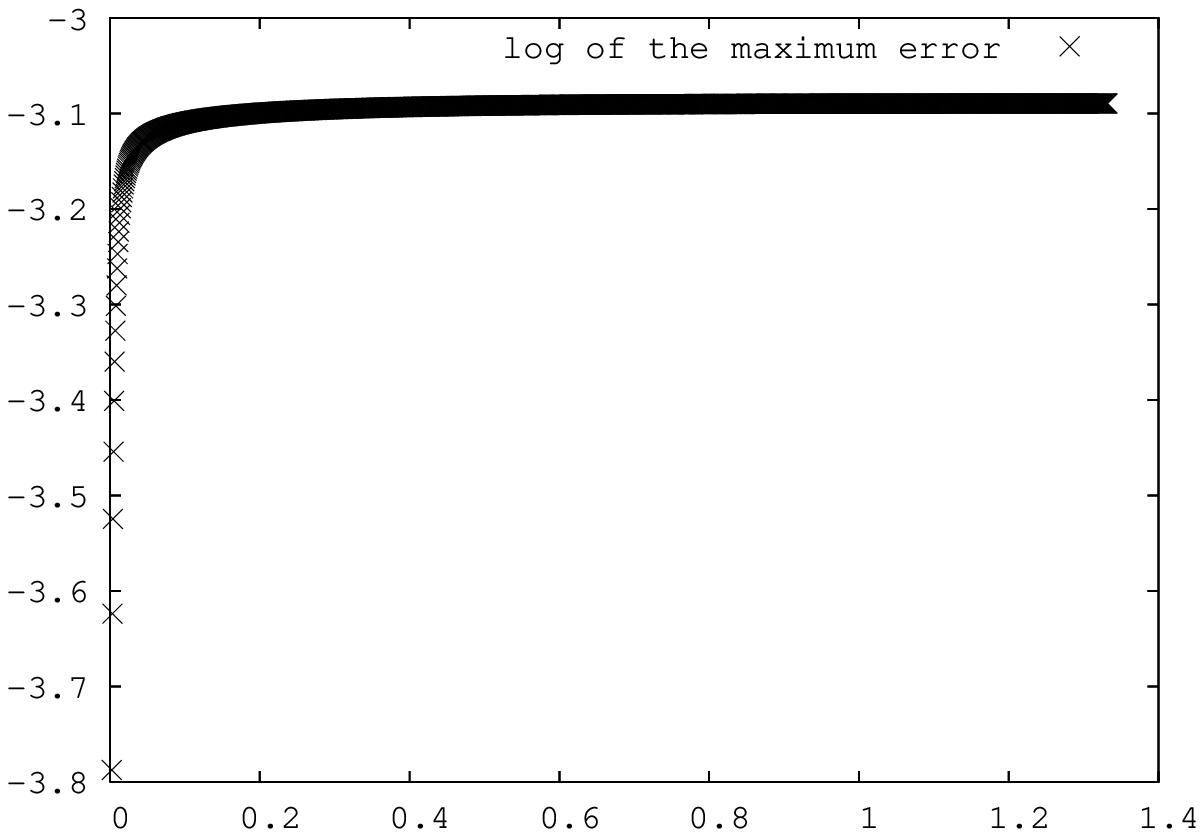}    }

\noindent  {\bf Figure 6}.  \quad 
Two-dimensional computation of the harmonic function  $ \, p_H(x,\, y)  = x^2 - y^2 \, $  
on a triangle with the D2T7 lattice Boltzmann scheme.  
No extra order is observed for the   $ {\rm L}^\infty $ error 
  when refining the mesh (left). 
The evolution in time is very  slow (right),  
even initializing the computation with  the exact solution !  
\smallskip \smallskip   

 \smallskip   \noindent 
In Figure~5, we have plotted the error field for the three versions 
(\ref{coefs-order2})(\ref{coefs-order4})(\ref{coefs-order6}) 
of the D2T7 lattice Boltzmann scheme. 
The results are qualitatively coherent: the more the scheme is theoretically precise, the
more the error is reduced. In Figure~6 (left), we observe that the  $ {\rm L}^\infty $ error 
is substantially reduced when the parameters induce a better precision. But the order of
convergence remains very close to second order even for  ``fourth order'' and 
``sixth order'' versions of the scheme using the set of parameters (\ref{coefs-order4})
or  (\ref{coefs-order6}).
In this case this default can be due to a possible deficit of time steps
and to crude boundary conditions. The lattice
Boltzmann scheme is explicit and the time iterations (see Figure~6, right) 
take too much time to reach convergence to the stationary state with a  
 satisfactory reduction of the error. 

 \smallskip  \monitem      Dissipation of a triangular Dirichlet mode

\noindent 
We have also experimented the relaxation of a Dirichlet mode. The first mode 
is simply a product of three ``sinus'' functions, as first explicited
by Lam\'e (see McCartin \cite{MC03}). 
With our nomenclature, the eigenvalue number  $ \, \ell \,$ 
is proportional to $ \, 3 \, (\ell-1)^2 \,$.  The reference value is 
in consequence equal to 12, 48 and 108 
for $ \ell $ equal to 3, 5 and 7 respectively. 
We present in Figure~7 the results at 
$ \, T = 4/3  \,$ and the evolution of the physical field at the center. The asymptotic
analysis obtained by successive mesh refinements is presented in Figure~8. 
We measure the error in time for the center vertex as the mesh size tends to zero and the 
 $ {\rm L}^\infty $ error at the precise  time $ \, T = 4/3  \,$  in the same conditions. 
The results are correct but not easy to interpret. The ``second order'' scheme is 
just a bit better that the order  $3/2.$ The ``fourth order'' version is of order~3
and the ``sixth order'' scheme hesitates 
between the orders~3 and~4.

\smallskip  

\noindent          
\includegraphics[width=.27 \textwidth] {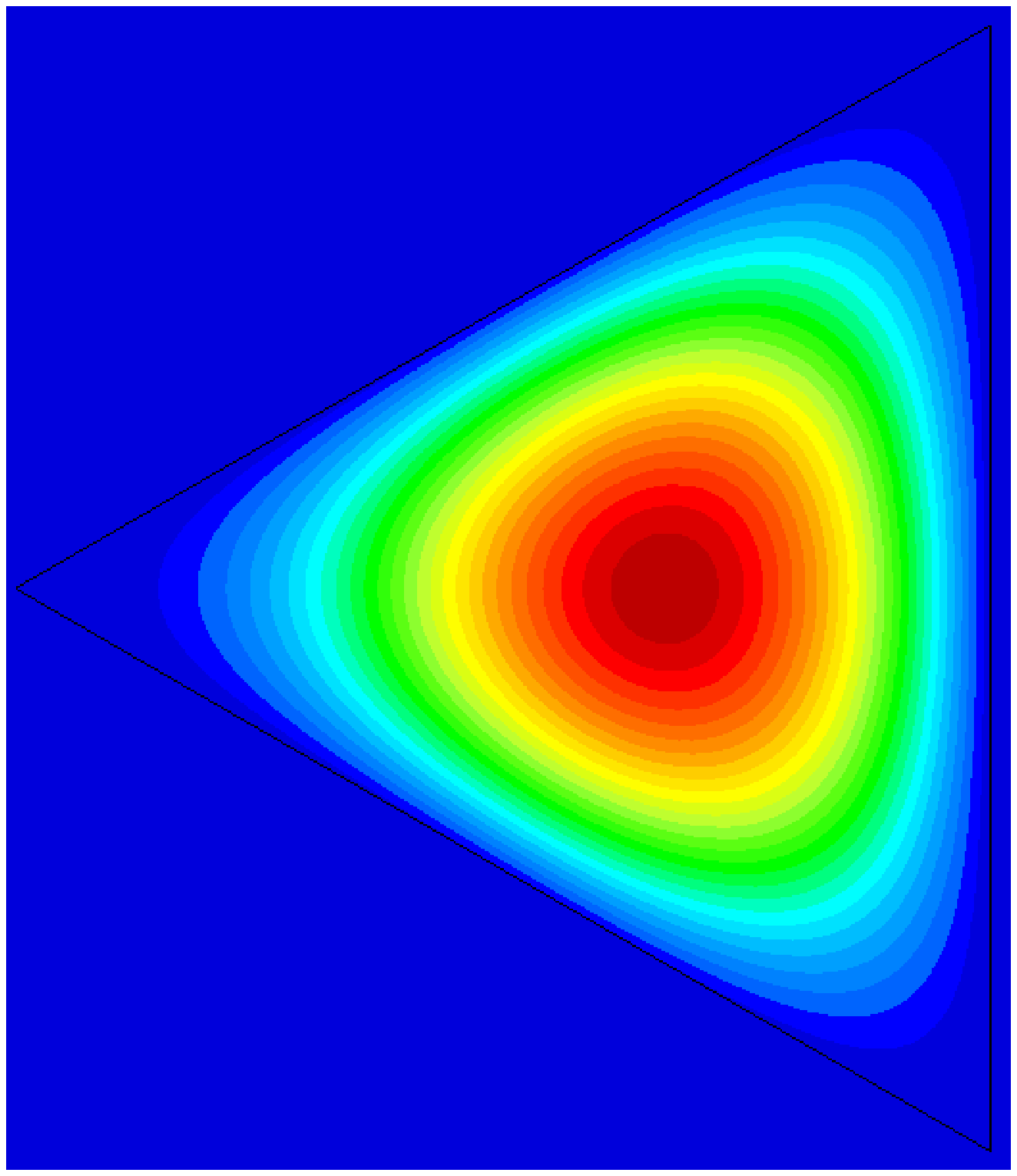} $\,\,$ 
\includegraphics[width=.27 \textwidth] {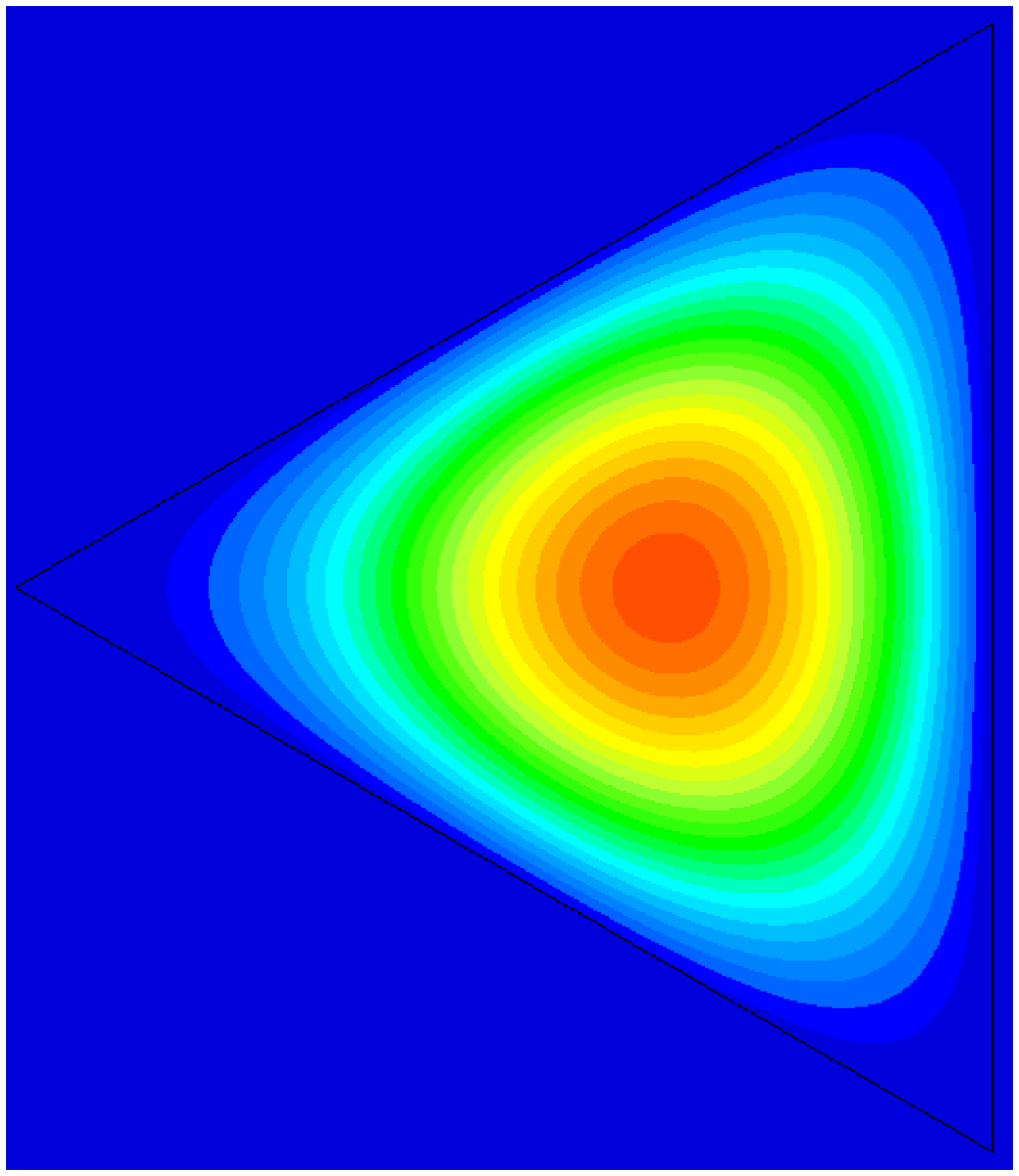}    $\,\,$ 
\includegraphics[width=.45 \textwidth] {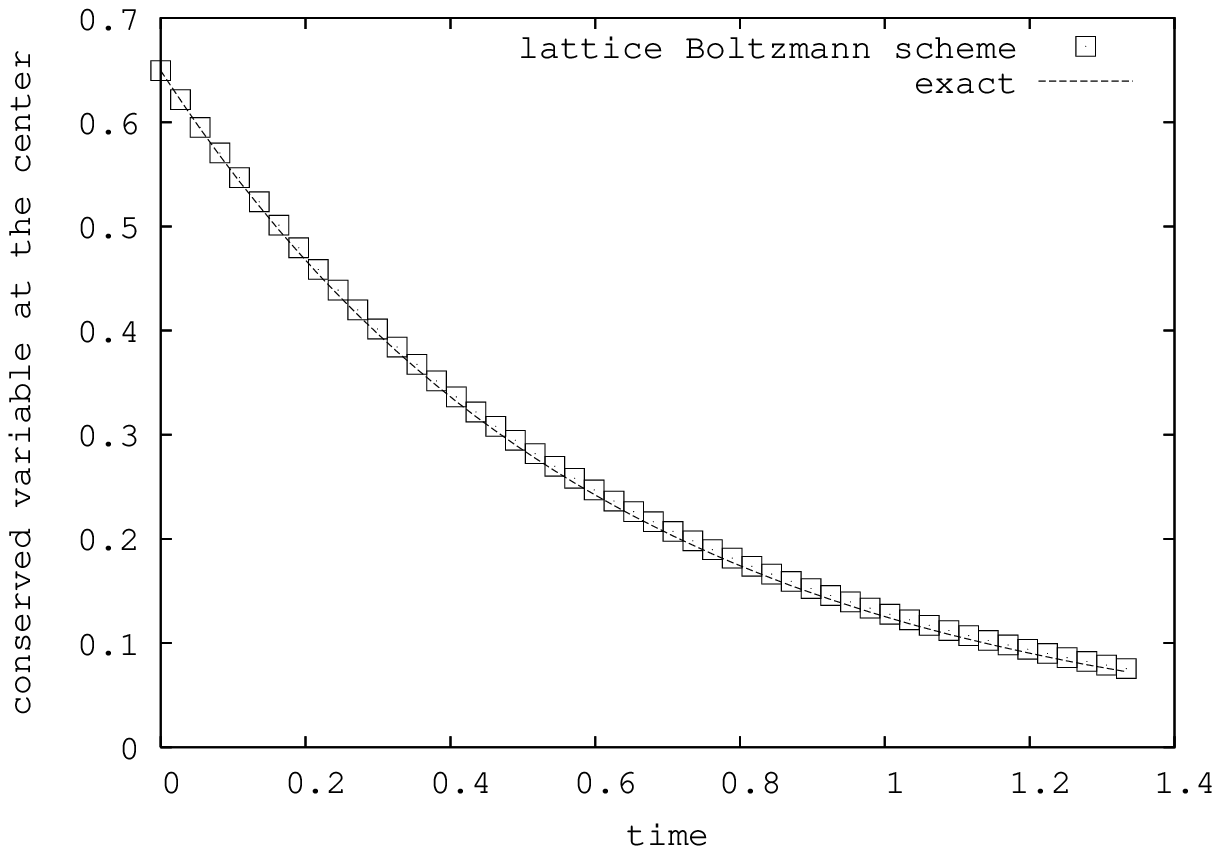}

\noindent  {\bf Figure 7}. \quad  
 Isovalues of the first Dirichlet mode for an equilateral triangle (left). 
  Dissipation of this  mode by  time evolution:  D2T7 solution at time $ \, T = 4/3 \,$ for 
76~points on the edge (2926~vertices, middle). 
Exponential decay at the center of the mesh (55~vertices, right).  
\smallskip \smallskip 

\smallskip             
\centerline { 
 \includegraphics[width=.35 \textwidth] {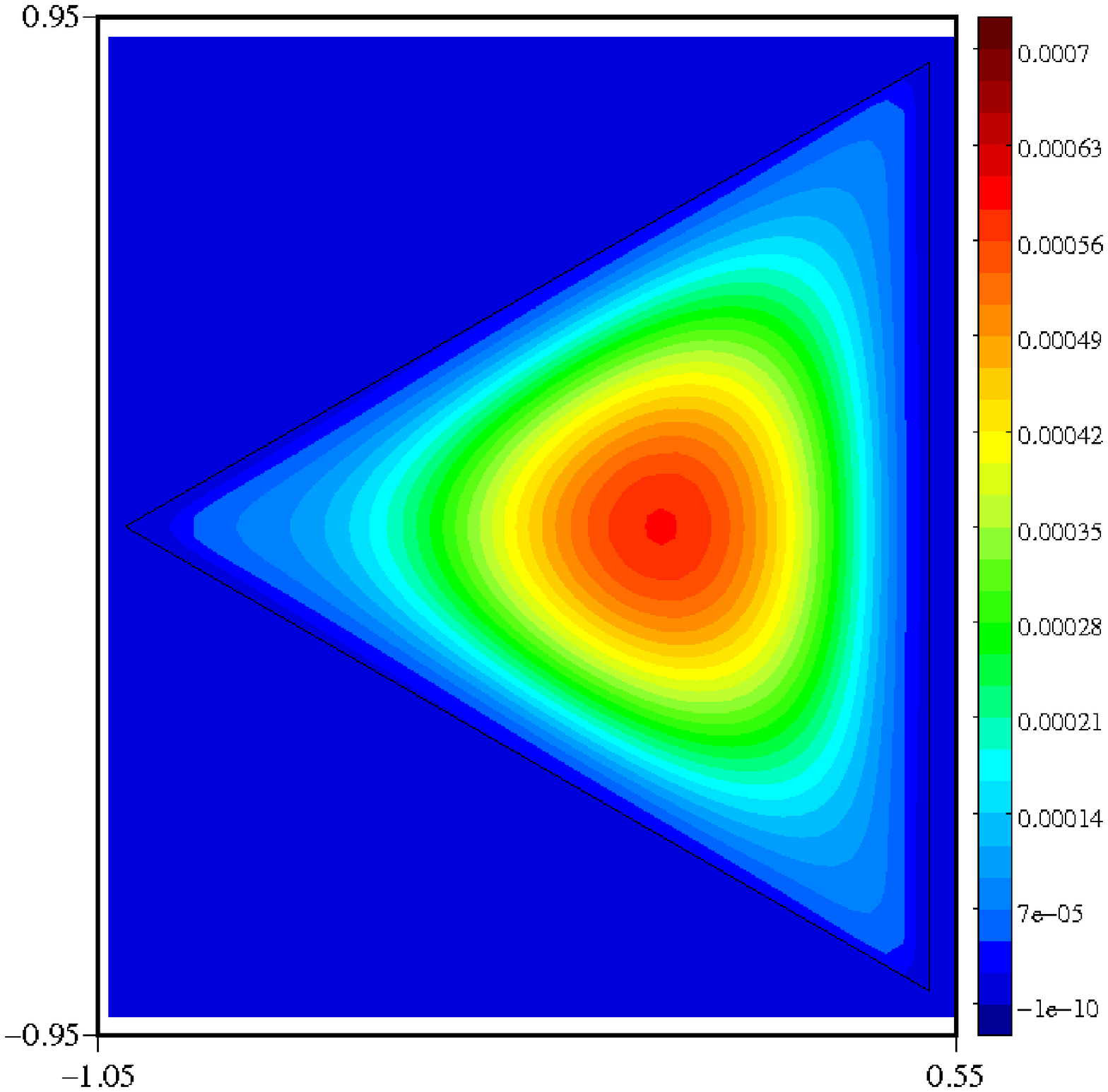} 
\qquad    
 \includegraphics[width=.49 \textwidth] {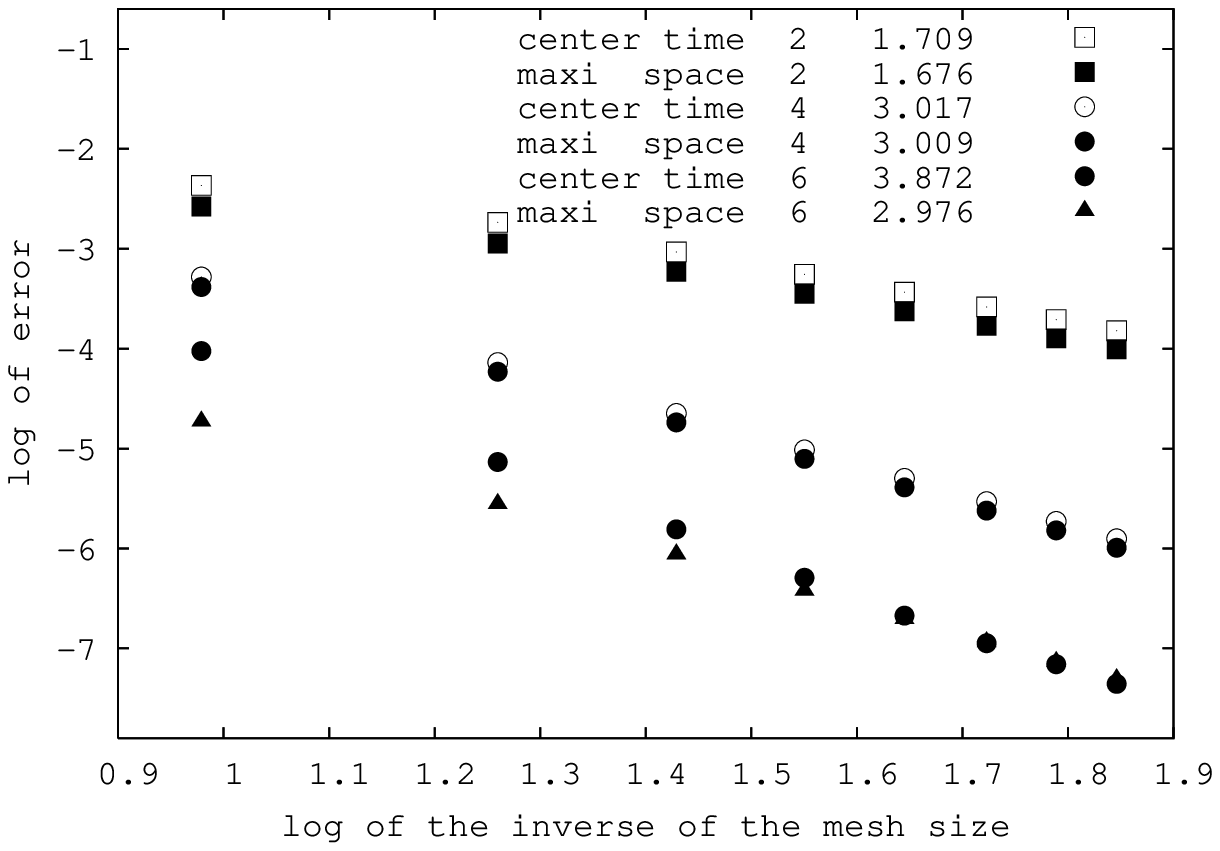}  }
 
\noindent  {\bf Figure 8}. \quad  
  Dissipation of the first Dirichlet mode. 
Isocontours of the 
field  of error  at $ \, T = 4/3 \,$ with a mesh composed by 61 points 
on the edge (left).   Negative values in blue and positive ones in red. 
 Time and $ {\rm L}^\infty $ space errors 
for several meshes and several ``orders'' with the D2T7 lattice Boltzmann scheme. 
The obtained accuracy is not 
the one proposed by the Taylor expansion method. The space numerical accuracy
is going from 1.7 to 3.0 with a good tuning of the numerical parameters.  
\smallskip \smallskip 

  \smallskip   \monitem      Dirichlet modes for a triangle 

\noindent 
We used the D2T7 lattice Boltzmann scheme (\ref{LB-scheme}) to define a linear operator 
$ \, f(t) \longmapsto f(t+\Delta t) \equiv A \smb f(t) \,$ 
where $\, f(t) \,$ is the vector of all unknowns for the entire mesh.
Then the first eigenvalues of the linear operator $\, A \,$ are computed with the Arnoldi
algorithm~\cite{Ar51}.

\smallskip                 
\centerline { \includegraphics[width=.35 \textwidth, angle=90] 
 {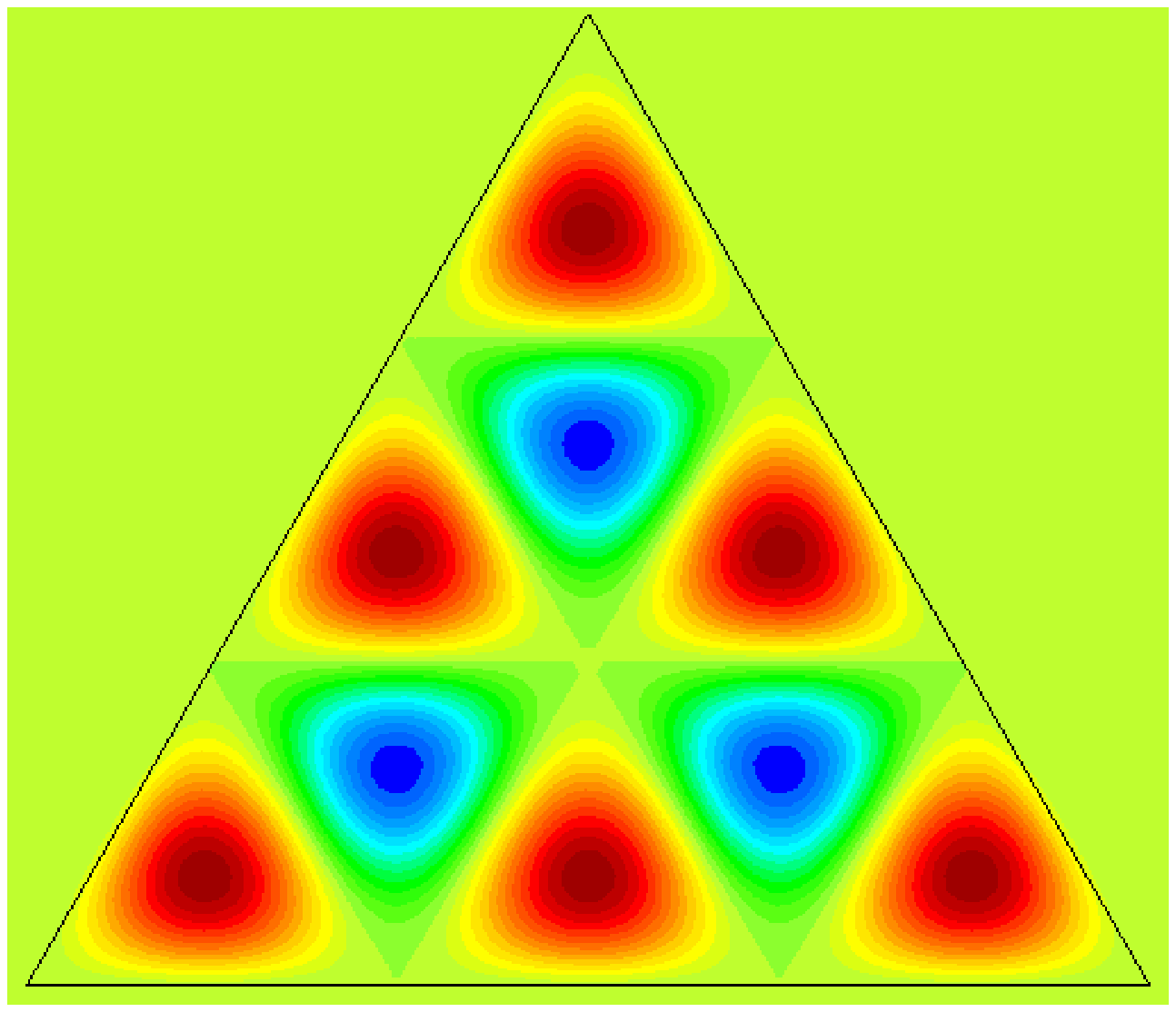} \quad 
 \includegraphics[width=.35 \textwidth, angle=90] {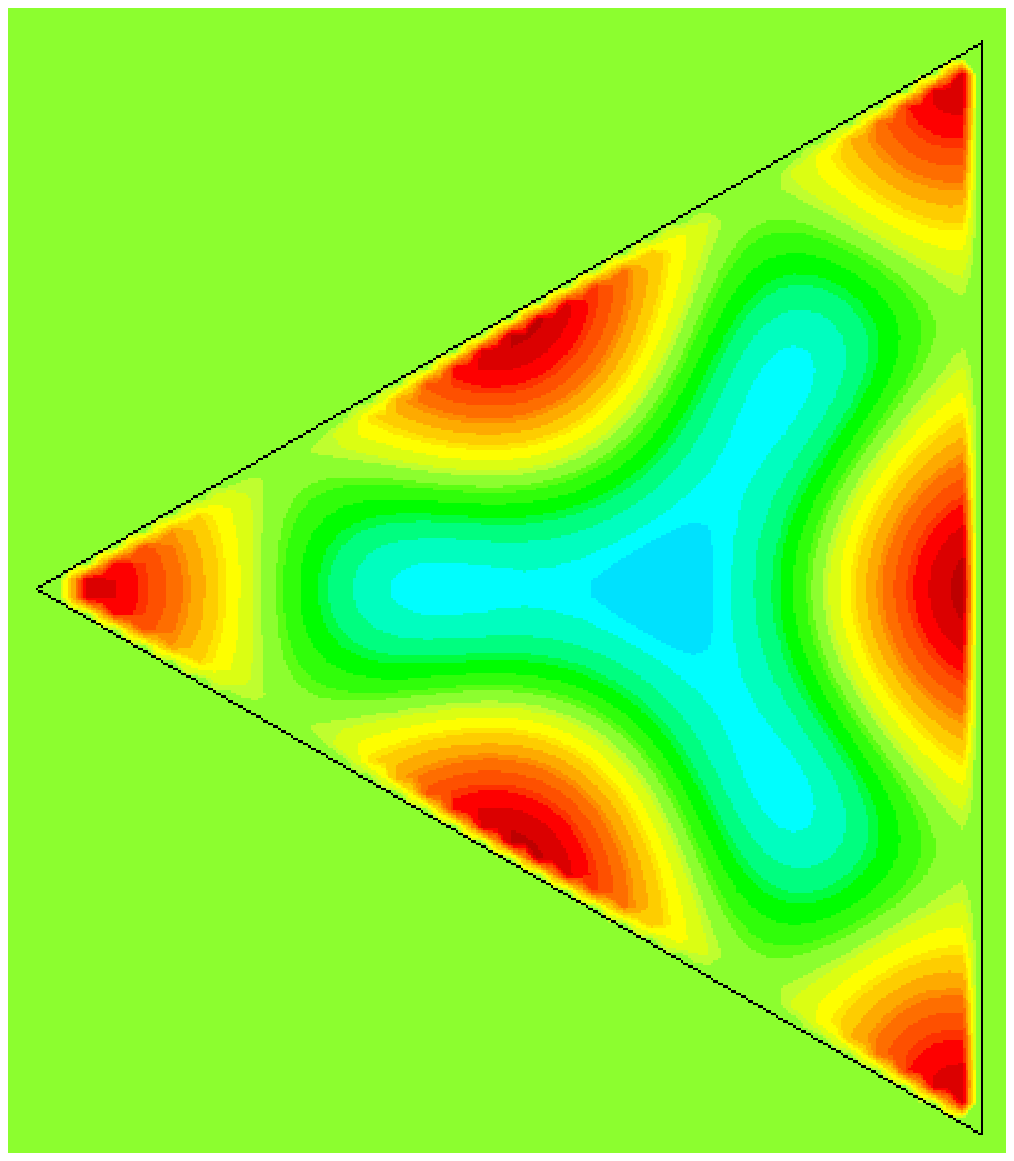} \quad 
 \includegraphics[width=.35 \textwidth, angle=90] {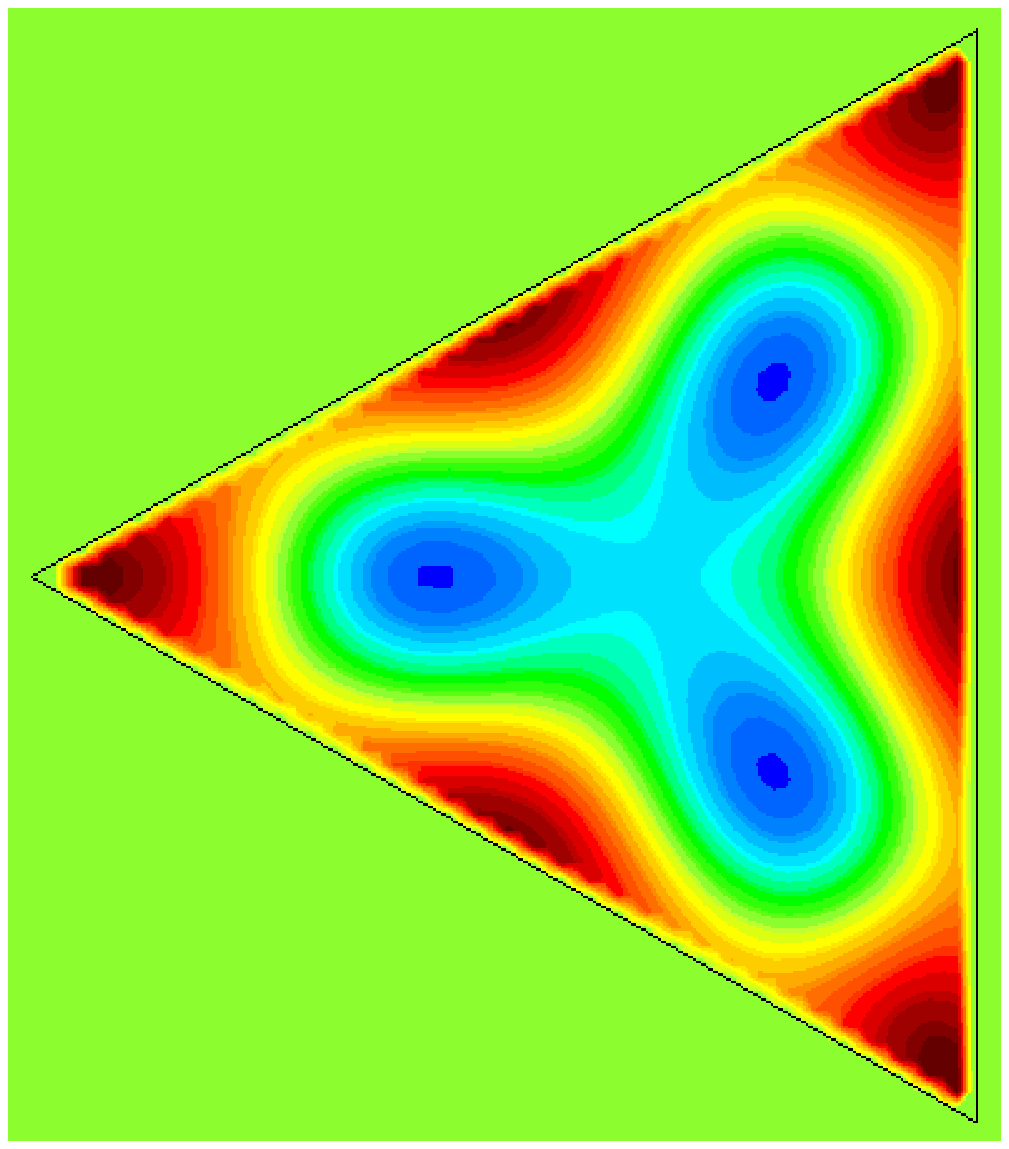}   }

\noindent  {\bf Figure 9}. \quad 
Isovalues of the 
Dirichlet mode  number  ``3''  of an equilateral triangle (left) and
errors for  a D2T7 computation. The exact reference eigenvalue 
is equal to 12 (in appropriate units). 
 The numerical eigenvalue is equal to  11.99902 
with second order parameters    (middle)
and to   11.99938  with ``fourth order'' parameters  (right). 
 The  $ {\rm L}^\infty $ error for the modes 
 is equal to  $ 4   \, 10^{-4}  $   at order 2 and 
 $    10^{-4}  $   at order 4. The figures show the isovalues of the error for both 
computations with different scales. We observe that the global shape of these 
errors is similar to isovalues of the reference eigenvector.
\smallskip \smallskip 
  
 \smallskip   \noindent 
Some exact reference modes are displayed in the left part of 
Figures~9 to~11. The numerical approximation
is globally of very good quality and we have plotted  the errors for different modes 
computed on the same lattice in the same figures.   
 We perform the computations for each mode, one with the
``second order'' accurate version of the D2T7 scheme and the other one with a 
``fourth order'' accurate tuning of numerical parameters. In each case, we compare 
the theoretical eigenvalue after applying a suitable normalization  
and the computed eigenvalue by the Arnoldi algorithm. 
The results are of good quality and the quartic parameters give a better precision 
for the numerical results. Even if the fourth order convergence is not established, 
the tuning  of parameters improves clearly the numerical quality.

\smallskip                  
\centerline { \includegraphics[width=.35 \textwidth, angle=90] 
{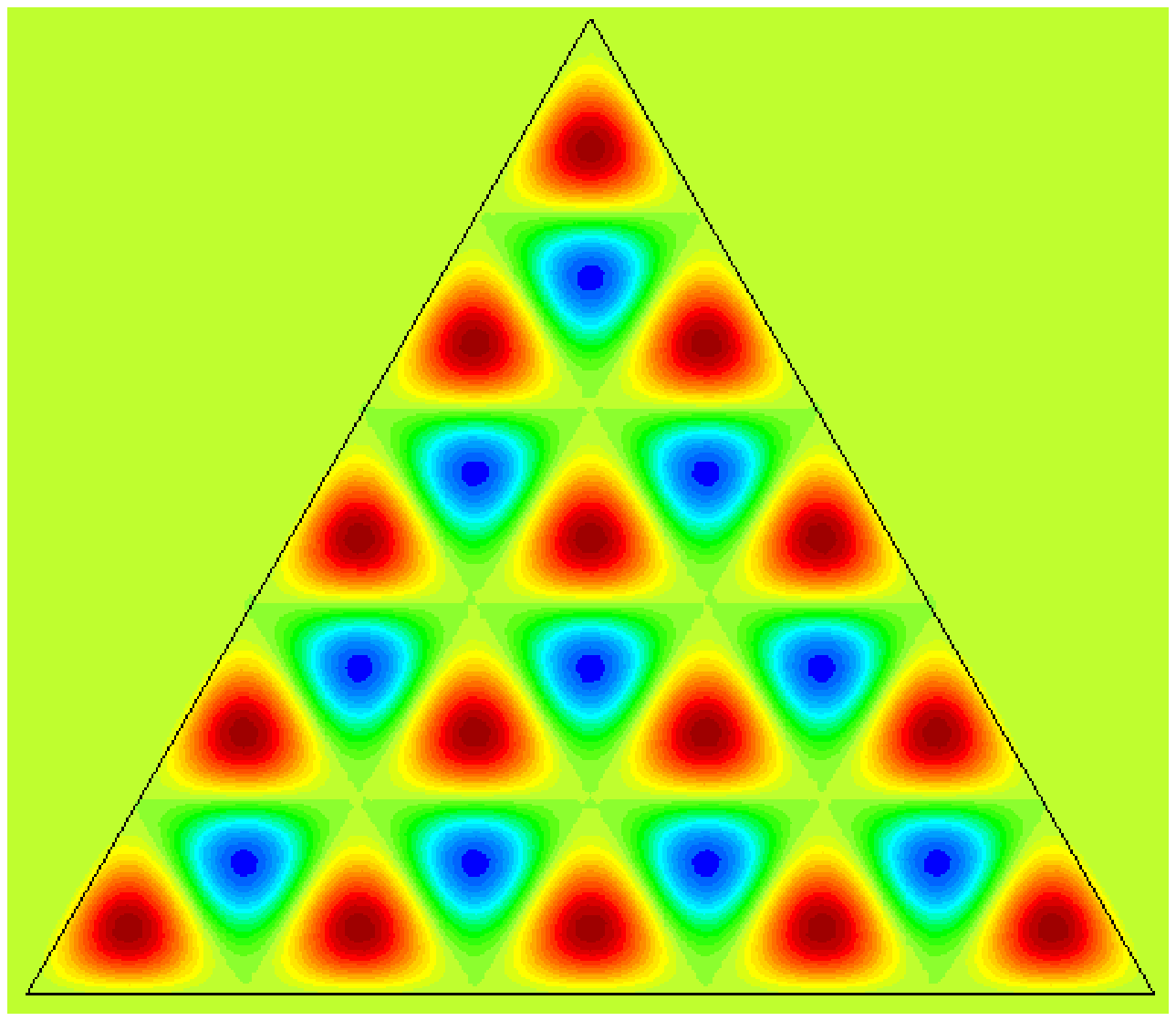} \quad 
 \includegraphics[width=.35 \textwidth, angle=90] {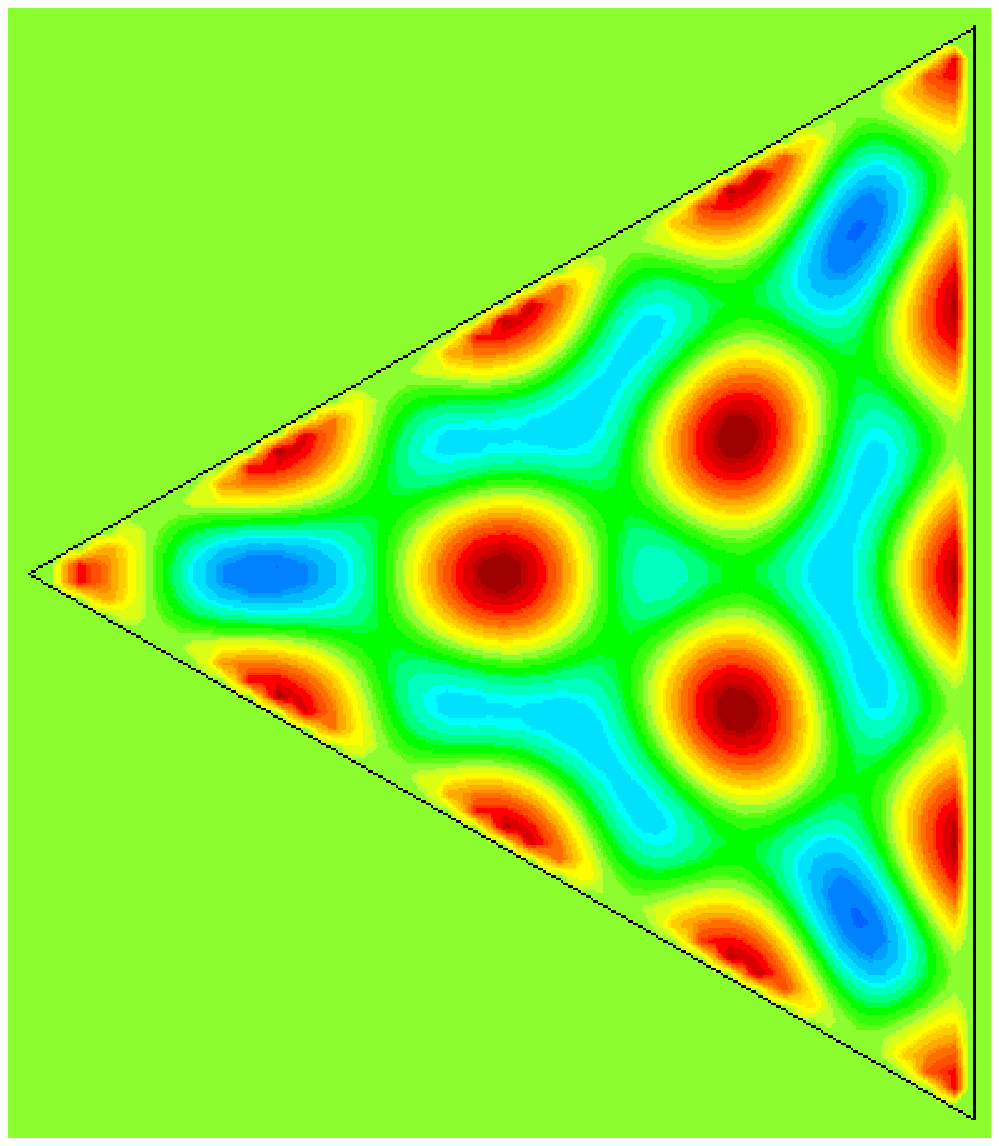} \quad 
 \includegraphics[width=.35 \textwidth, angle=90] {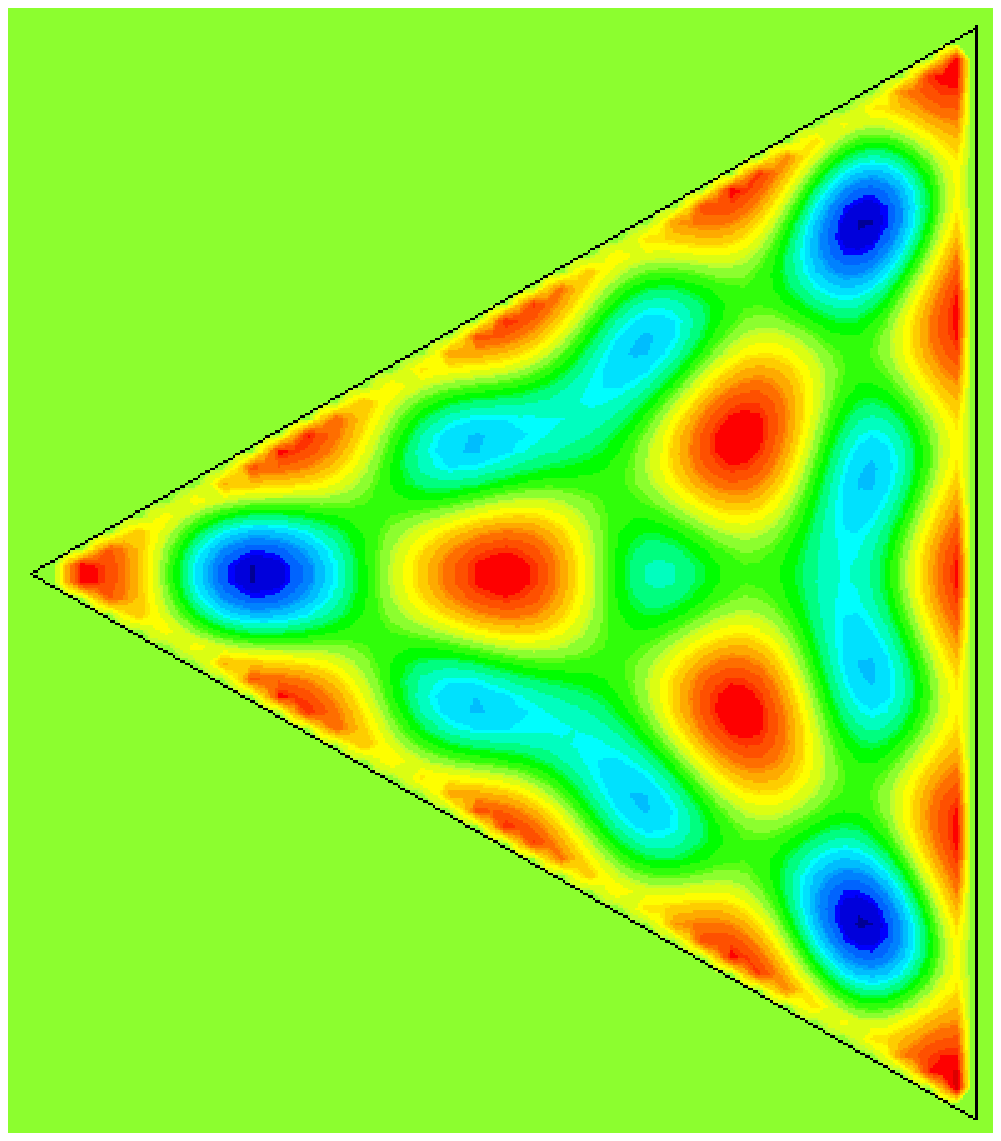}   }

 \noindent  {\bf Figure 10}. \quad 
Isovalues of the 
Dirichlet mode  number  ``5''  of an equilateral triangle (left) and
errors for  a D2T7 computation. The exact reference eigenvalue 
is equal to 48. The numerical eigenvalue is equal to   
47.98339 at order 2  (middle) 
and to   47.98842  at order~4  (right).   The  $ {\rm L}^\infty $ error for the modes
 is equal to    $ \, 3 . \,  10^{-2} \,  $  at order 2 and 
 $ \,  1.1   \,  10^{-3} \, $ at order~4.  
\smallskip \smallskip 

\smallskip                   
\centerline { 
\includegraphics[width=.35 \textwidth, angle=90] 
 {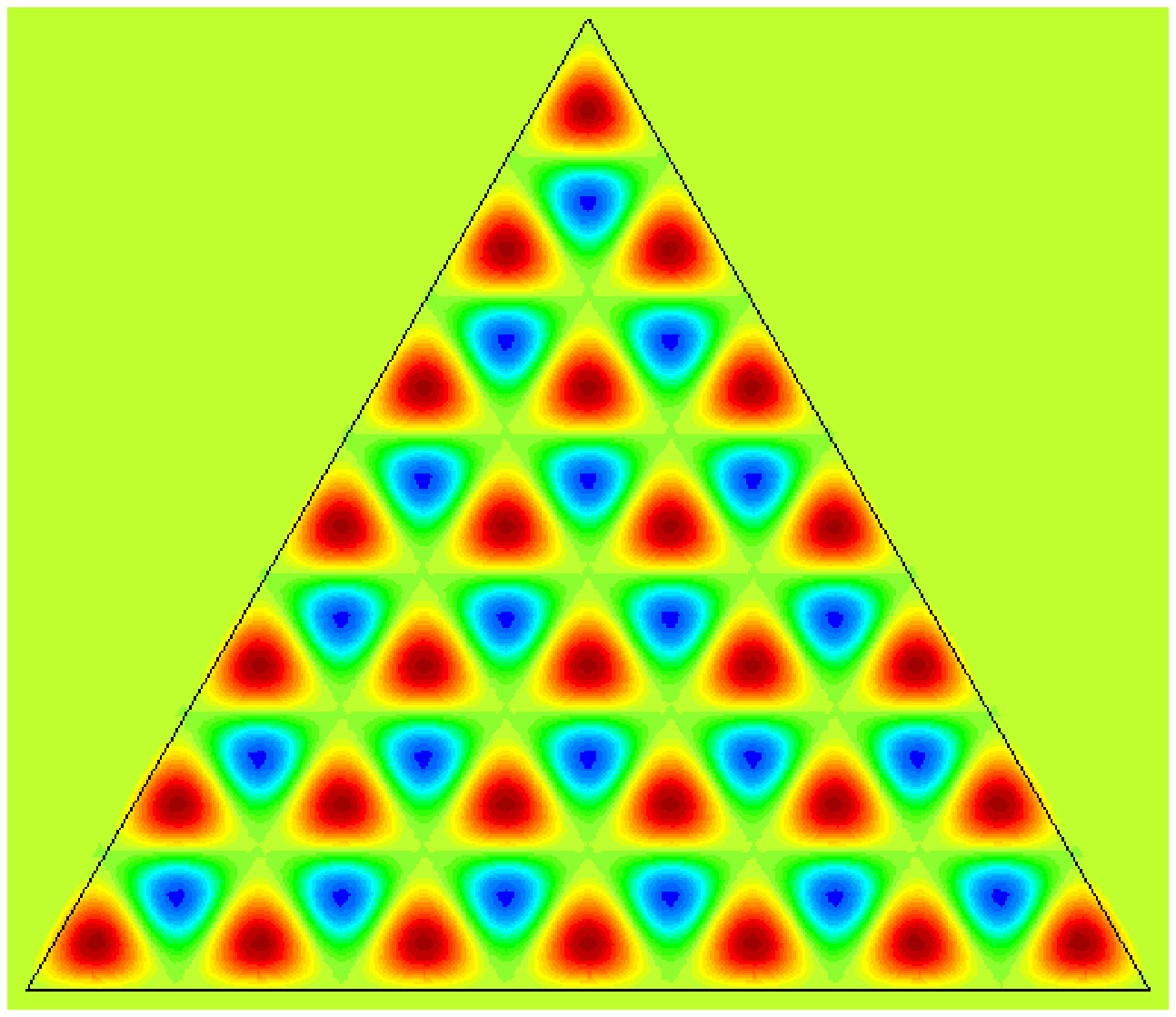} \quad 
 \includegraphics[width=.35 \textwidth, angle=90] {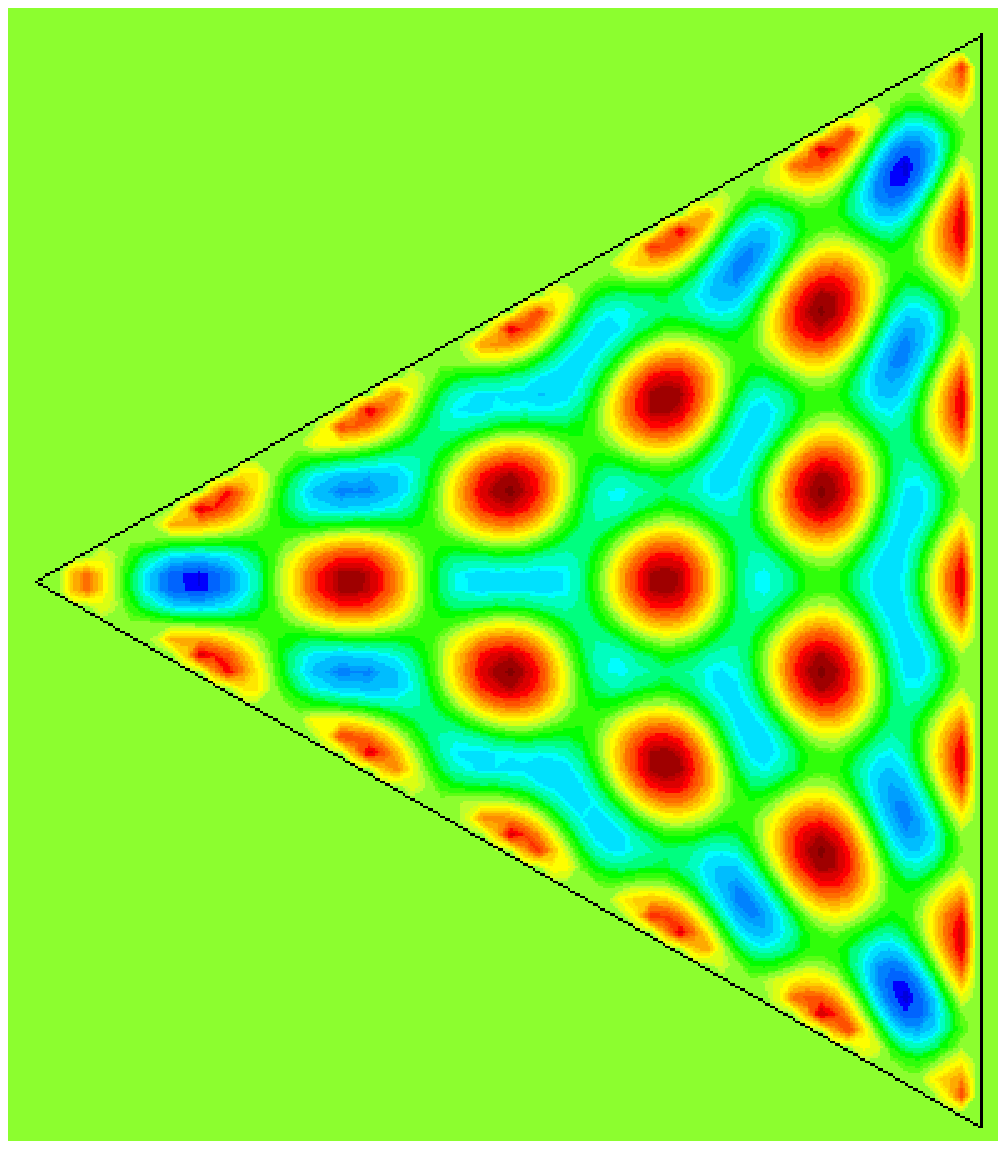}  \quad 
 \includegraphics[width=.35 \textwidth, angle=90] {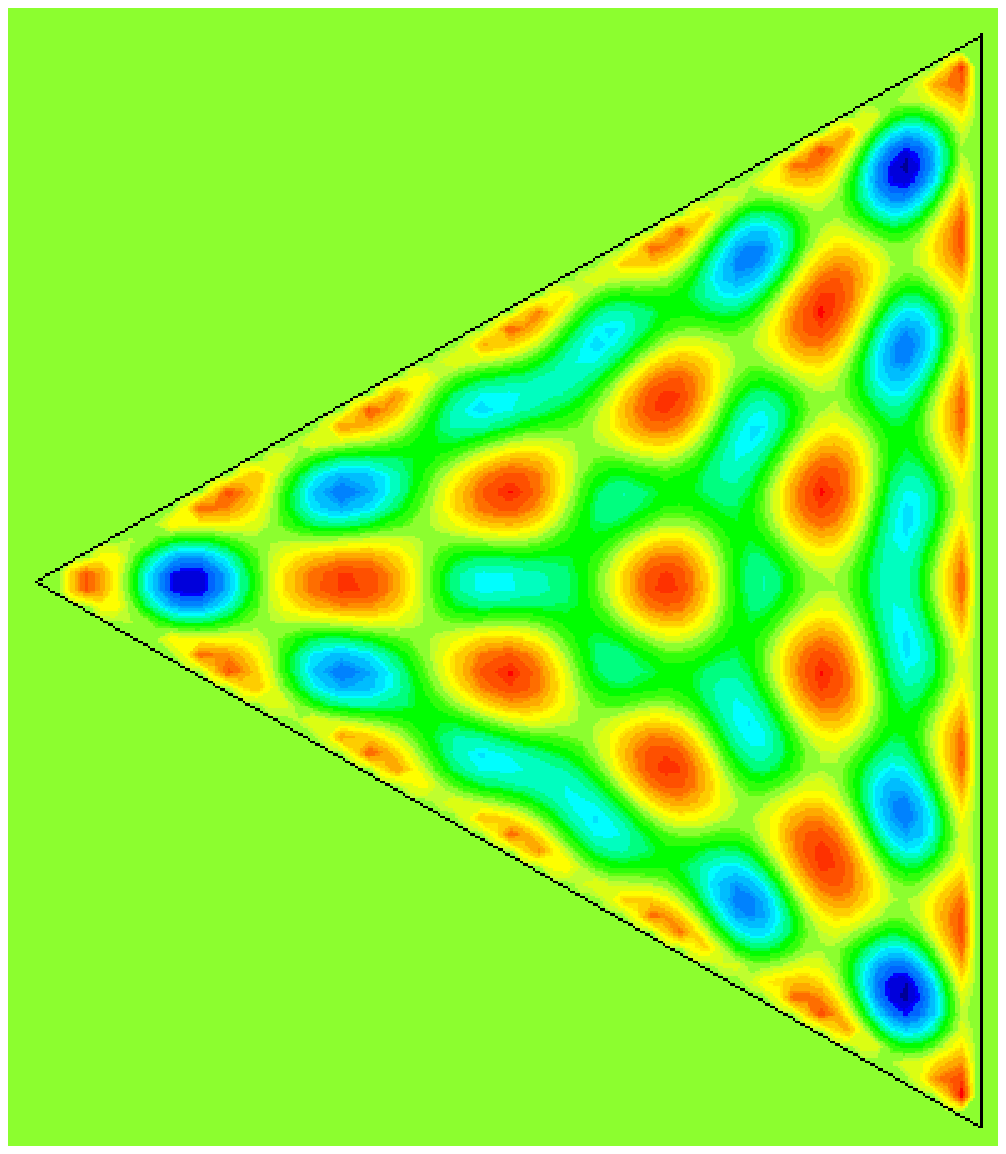}   }

\noindent  {\bf Figure 11}. \quad 
Dirichlet mode  number  ``7''  of an equilateral triangle (left) and
errors for  a D2T7 computation. The exact reference eigenvalue  
is equal to 108. The numerical eigenvalue is equal to 107.90777 at order 2  (middle)
and to  107.92705  at order 4  (right).  The  $ {\rm L}^\infty $ error for the modes with 
 this computation is equal to   $ \, 1.02  \,\,  10^{-2}  \, $   at order 2 and 
 $ \,  4.2   \,  10^{-3} \, $  at order~4.   
\smallskip \smallskip 

  \newpage 
\bigskip \bigskip   \noindent {\bf \large 5) \quad Lattice Boltzmann scheme on arbitrary  meshes ?}  


  \noindent   
Imagine that we move the vertices in the Bravais lattice  presented at Figure~4 (left). 
We obtain a topologically regular  mesh 
in the sense that the number of edges containing a given internal vertex is constant. 
An example is  proposed at Figure~12.  
%
This mesh is a good candidate for future extensions of the lattice Boltzmann scheme. 
With this kind of classical finite element type mesh, it is possible
to use all the engineering    
tools of automatic meshing in two and three space dimensions
as described {\it e.g.} in \cite {Ge92}. 
But this goal is still not the purpose of the present contribution. 
The vertices of the mesh of Figure~12 
are now the nodes of a cellular complex 
and  each vertex has a {\bf constant} number of neighbours.
In other terms,   the degree of each vertex is constant.  
We denote by $ \,  x_j \equiv x + \xi_j(x) \, \Delta x \, $ the  vertex belonging 
to the lattice $ {\cal L}$ with a local neighbouring  number   $ \, j \, $ 
 relative  to  the  vertex  $ \, x$.
Remark that the vertex   $x$ is also a  neighbour of the vertex $x_j $ 
with a local number   $ \ell \equiv n_j(x) $. 
We have the obvious relation  $ \, \xi_j(x) + \xi_\ell (x_j) \equiv 0 \,$ 
and in other terms the identity
\moneq   \label{local-dual} 
\xi_j(x) + \xi_{n_j(x)}  (x_j)   \equiv 0  \, . 
\monend   
As previously, we denote by  $  \, f_j (x) \, $ the density of particles going from 
  vertex    $\, x \, $   towards   vertex $ \, x_j $.   
Moreover, the   outgoing   particles from   vertex  $ \, x_j \, $  
are also  ingoing   particles  ``into''   vertex   $ \, x \, $   
with an index denoted 
by $ \, \ell$.

\smallskip                    
\centerline { \includegraphics[width=.33 \textwidth, height=.33 \textwidth]
 {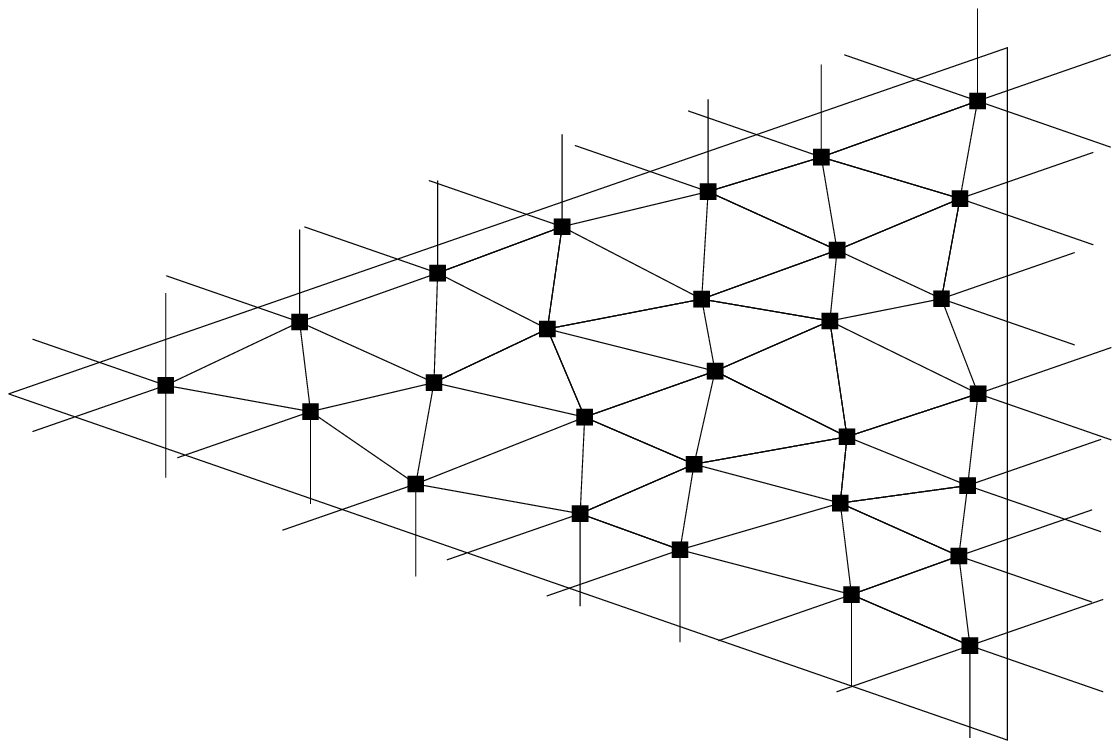} }  

\noindent  {\bf Figure 12}. \quad
Triangular lattice obtained from a little random displacement 
of the vertices of an equilateral triangular mesh.  
\smallskip \smallskip  
 
 \smallskip    \noindent  
We precise the previous notation. 
If   $   \, f_j (x) \, $ is the  density of particles from 
the vertex    $\, x \, $  towards the vertex $\, x_j , $ 
$  \,  f_j^* (x) $ denotes the same quantity  {\bf after relaxation}.
In a dual vision, we denote by    
$  \,  f_\ell^* (x_j) \, $ the density of particles going from the vertex
   $ \, x_j \, $ in the direction of the  vertex   $\, x \, $
 after relaxation. 
We have also to consider the  density $  \, \widetilde {f_j}(x) \, $
of particles going from the vertex $x_j $ 
 towards the vertex   $ \, x$. 
The lattice Boltzmann scheme is a particle method.
The   flight of   particles between the vertex   $ \, x_j \, $ 
and the vertex      $\, x \, $   takes exactly one time step :  
   $ \,\, \widetilde{f_j}(x,  \, t+\Delta t) \,=\,  f_\ell^* (x_j,  \, t)  $.  
If we replace the notation $ \, \ell \,$ for the index of vertex $x$
relative to its neighbour $ \, x_j \,$ 
by the notation $ \, n_j(x) \,$ introduced previously at relation 
(\ref{local-dual}), the lattice Boltzmann scheme takes the form
\moneq   \label{LB-particle} 
 \widetilde{f_j}(x, t+\Delta t) \,=\,  f_{n_j(x)}^* (x_j, t) \, .  
\monend   
In the case of general meshes, the relation (\ref{LB-particle})
replaces the initial formula (\ref{LB-scheme}), 
correct only for Bravais lattices, as illustrated in Figure~13.

\smallskip                    
\centerline { \includegraphics[width=.33 \textwidth] {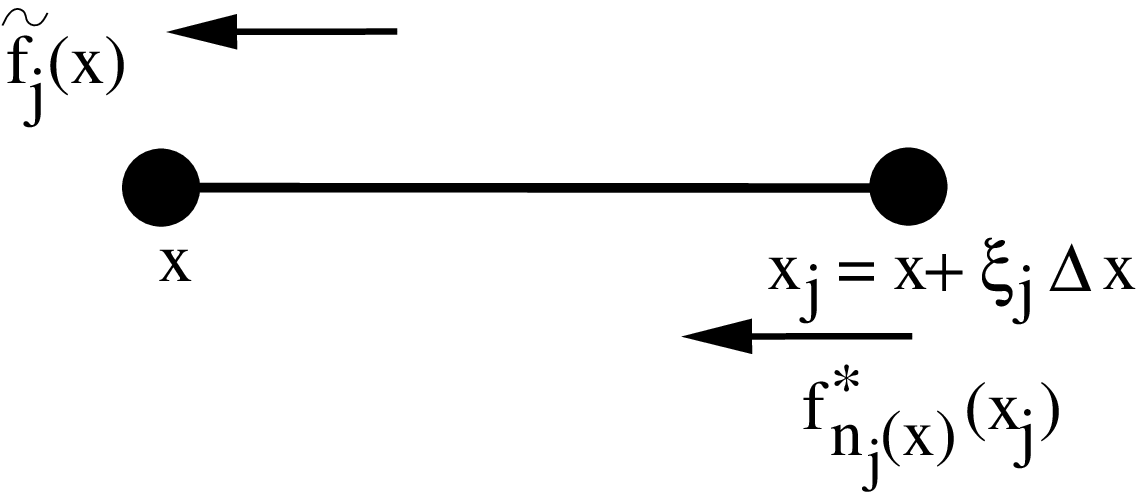} }    

\noindent  {\bf Figure 13}. \quad
Iteration of a lattice Boltzmann scheme: the
ingoing particles $\, \widetilde {f_j}(x) \,$ into  vertex $x$ are coming from 
the neighbouring vertex $ \, x_j   \,$ after a relaxation step.  
\smallskip \smallskip 
 
 \smallskip    \noindent 
We precise now how to compare the ingoing particles 
$ \,  \widetilde{f_j}(x) ,\,$ the outgoing particles 
$ \,  f_j^* (x) \,$ emitted from the vertex $ \, x \, $ 
and the associated moments. 
We introduce first a matrix $ \,  \widetilde{M}(x) \,$ 
in order to compute the moments $ \, m_k(x) \,$ from 
the ingoing particles. As previously (see the relation 
(\ref{polynomes})), we suppose given a family $ \,  {\cal P} \,$ 
of polynomials $ \, p_k .\,$  In an analogous way suggested 
by the relation (\ref{matrice-M}), we just  reverse the direction
of velocities  
and we have
$\,\, \widetilde{M}(x)_{k j} \,=\, p_k \big( - \xi_j \big) ,\,$ 
with $ \,  p_k \in  {\cal P}  .  \,$  
If the polynomials $ \, 1 ,\, $  $  X \, $ and   $ \,  Y \, $ are the first
polynomials of the family  $ \,  {\cal P} ,\,$ we have as 
in the previous studies 
$ \,  \, 
\widetilde{M}(x)_{\, 0  j} = 1 \,, \,\,  
\widetilde{M}(x)_{\, \alpha j} =  -\xi_j^\alpha (x) \,,
\,\,  1 \leq \alpha \leq d . \, $ 
The moments are evaluated for the incoming particles  with the natural relation
\moneq   \label{moments-entrants}  
m_k(x) \equiv \sum_j  \widetilde{M}(x)_{k j} \, \, \widetilde{f_j}(x) \,,
\quad 0 \leq k \leq q-1 \,, \quad x \in  {\cal L} \, .   
\monend   
The relaxation step is essentially unchanged. The moments can be seen 
as the eigenvectors  of the jacobian of the relaxation matrix 
(see {\it e.g.} \cite {Du10})  and this operator is diagonal with this representation: 
\moneq   \label{relax-general}  
m_k^* (x) = m_k (x) + s_k (x) \, \big( m_k^{\rm eq}(x)  - m_k(x) \big)  \, , 
\monend   
where the index $k$   in relation  (\ref{relax-general}) 
is running on all nonconserved moments. 
The outgoing particles after relaxation are supposed to 
be a {\bf linear} functional of the moments: 
\moneq   \label{f-star}   
f_j^*(x)  \equiv \sum_k   P(x)_{j k}  \,\,  m_k^*(x)  \, , 
\quad 0 \leq j \leq q-1 \,, \quad x \in  {\cal L} \, .   
\monend   
The question is now to determine the matrix $P$.  
We have the following property.

\smallskip  \noindent  {\bf Proposition 1. \quad Transition matrix from moments to particle
distribution.}

\noindent  
 If the Taylor expansion approach is valid at the order zero and if 
each internal node of lattice  $\,  {\cal L} \, $ is of constant degree (the number 
of neighbours of each vertex does not depend of the vertex $\, x \in {\cal L}$), the
matrix $\, P(x) \,$   of relation (\ref{f-star}) is given by the relation 
\moneq   \label{matrice-P} 
P(x)_{\, i \,  \ell}   \,=\, 
  \big(  \widetilde{M}(x_i) \big)^{-1} _{\,\,\,   n_i(x) \, \ell}  \, . 
\monend   

\noindent  {\bf Proof of Proposition 1}.

 \noindent
The proof can be conducted as follows. We start from the time iteration
(\ref{LB-scheme}) of the lattice Boltzmann scheme. 
  Then  after multiplication by the matrix 
$ \,  \widetilde{M}(x) , $ 
with the help of  (\ref{moments-entrants}),  (\ref{LB-particle}) and 
 (\ref{f-star}), we have 

\smallskip \noindent 
$ \displaystyle \smash {  m_k(x, \, t + \Delta t ) =    \sum_j   \widetilde{M}(x)_{k j}  \,\,  
\widetilde{f_j}(x, \, t + \Delta t) } $ 
$ \,  =   \displaystyle \,   \sum_j   \widetilde{M}(x)_{k j}  \,\,   f_{n_j(x)}^* (x_j,
\,t)  $

\smallskip \noindent 
$     \displaystyle  =   \, \smash {   \sum_j   \widetilde{M}(x)_{k j} 
\,  \sum_\ell   P(x_j)_{n_j(x) \,  \ell}  \, \, m_\ell^* (x_j,  \,t)  } $  
$   \displaystyle     \,  =   \,   \sum_\ell \, \big( \sum_j   \widetilde{M}(x)_{k j}  \, 
 P(x_j)_{n_j(x) \,  \ell}    \big) \,\,  m_\ell^* (x + \xi_j \Delta x , \, t) . \,  $ 

\smallskip \noindent 
We expand this relation at order one. Due to relaxation, we just have
a small perturbation between $m$ and $m^*$ : 
\moneqstar 
 m_k(x) + {\rm O}(\Delta t) =  m_k^*(x)  + {\rm O}(\Delta x) .  
\monendstar 
In consequence,   
\moneqstar 
   \sum_j   \widetilde{M}(x)_{k j}  \, \, 
 P(x_j)_{n_j(x) \,  \ell}   \, \equiv  \, \delta_{k \, \ell}  
\monendstar 
and in other terms, 
\moneqstar   P(x_j)_{\, n_j(x) \,  \ell}  =  
  \big(  \widetilde{M}(x) \big)^{-1} _{\,\,\, j\, \ell}   . 
\monendstar  
We change the names  of the vertices. 
We replace the letter $ \, x_j \, $ by the letter $\, x$. 
Then we  replace the index $ \, n_j(x) \, $ by some neighbor $\, i \, $ of vertex
 $\, x \, $ and the index $ \, j \,$ is now equal to $ \, n_i(x) .\, $  
With this change of notation,  we obtain      $ \, \,      P(x)_{\,i \,  \ell}  =
  \big(  \widetilde{M}(x_i) \big)^{-1} _{\,\,\,   n_i(x) \, \ell}  \,\, $ 
which is exactly the relation~(\ref{matrice-P}).  \hfill $\square$

\smallskip \noindent 
We can now make explicit the   d'Humi\`eres  lattice Boltzmann scheme 
on an arbitrary    mesh where the degree of each vertex is constant.  
When all the outgoing densities of particles 
$ \,  f_j^* (x,\, t )  \, $ are known for all the vertices of the lattice
at some discrete time $ \, t ,\, $ the ingoing densities 
$ \,   \widetilde{f_j}(x,  \, t+\Delta t)  \, $ at the new time step 
are simply evaluated by a free flight (\ref{LB-particle}) during one  time step. 
Then the moments  $ \, m_k \,$ are a local linear transform 
of the particle densities thanks to (\ref{moments-entrants}). 
The first moments compose a set $ \, W(x) \, $ of conserved variables 
and the equilibrium moments  $ \,  m^{\rm eq} \,$  are a given (in general nonlinear)
function $ \, G(W) \,$ of this field: 
$ \,\,  m_k^{\rm eq}(x) = G_k ( W(x))   ,\, \,\,  x \in {\cal L} . \, $ 
The relaxation of moments follow the relation
(\ref{relax-general}). Note that in general the coefficients 
$ \, s_k(x) \,$ now depend {\it a priori} explicitly on the vertex $x$. 
Last but not least, the  outgoing particles at the new time
step   from the vertex    $\, x \, $ follow the local linear transform 
(\ref{f-star}).

\bigskip \bigskip   \noindent {\bf \large 6) \quad D2T4 scheme for equilateral triangles}  


\noindent   
We consider a general two-dimensional mesh  $ \, {\cal L} \, $  composed by triangles. 
Note here that    a cellular complex 
is   composed by ``vertices'' in $ \,   {\cal L}^0  \,$ of dimension zero,  
 by edges  in $ \,   {\cal L}^1  \,$ of dimension one and by 
triangles of dimension two:   $\, x \in  {\cal L}^2  . \,$ 
 In  other words, we adopt a ``cell center''
framework in the sense  proposed in Roache \cite{Ro72}. 
We can also locate the degree of freedom $ \, x \, $ 
at the center of gravity of the corresponding triangle. 
Remark that we make here {\it a priori} no other regularity hypothesis. 
Each triangle $\, x \, $ has three edges. Each edge 
inside  the border of  $ \, x \, $ is  part of  the boundary of (at most) two triangles :  
the triangle $ \, x \, $  itself  and its $\, j^{\rm th}\, $ neighbor $\, x_j $.   
It is then natural to consider outgoing particles 
$\, (f_j)_{0 \leq j \leq 4} \,$ going from $ \, x \, $ towards 
$ \, x_j \, $ with a local velocity 
$ \, \xi_j(x) \, \Delta t  \,$ chosen in such a way that the centers of both triangles
$ \,x \, $ and $ \, x_j \, $ are joined in exactly one time step of duration $\, \Delta t$. 
Of course, the null  velocity is not excluded. 
This remark explains  the name ``D2T4'' of this type of lattice Boltzmann scheme. %
A typical regular mesh for a D2T4 computation is presented in Figure~14.

\smallskip                    
\centerline { \includegraphics [width=.32 \textwidth, height=.33 \textwidth] 
{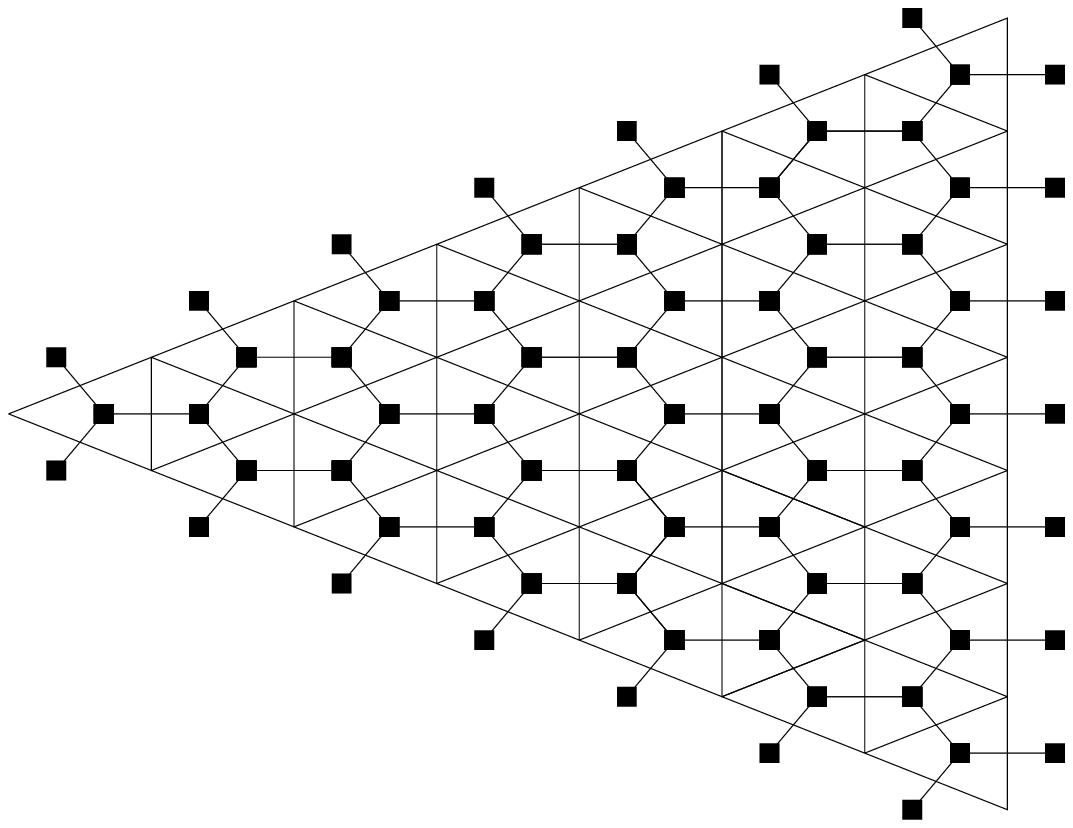} } 

\noindent  {\bf Figure 14}. \quad 
 Typical mesh with equilateral  triangles.  The four degrees of freedom of D2T4 scheme
 are  located at the center of  gravity of each triangle.   The links between triangles 
create the  dual hexagonal mesh around the vertices of the triangular mesh.  
\smallskip \smallskip 
 
  \smallskip   \noindent 
The degrees of freedom in  Figure~14 are the centers of the initial triangular mesh.
This ``secondary mesh'' is no longer  a Bravais lattice. We lose the possibility of
straight propagation of particles in the lattice and also  the symmetry property
of Bravais meshes emphasized in Figure~1. But we keep the property 
that the number of neighbours  
is constant. And this property is maintained
whatever the initial triangulation with cellular complexes. 

\smallskip \smallskip                    
\centerline { \includegraphics[width=.35 \textwidth] {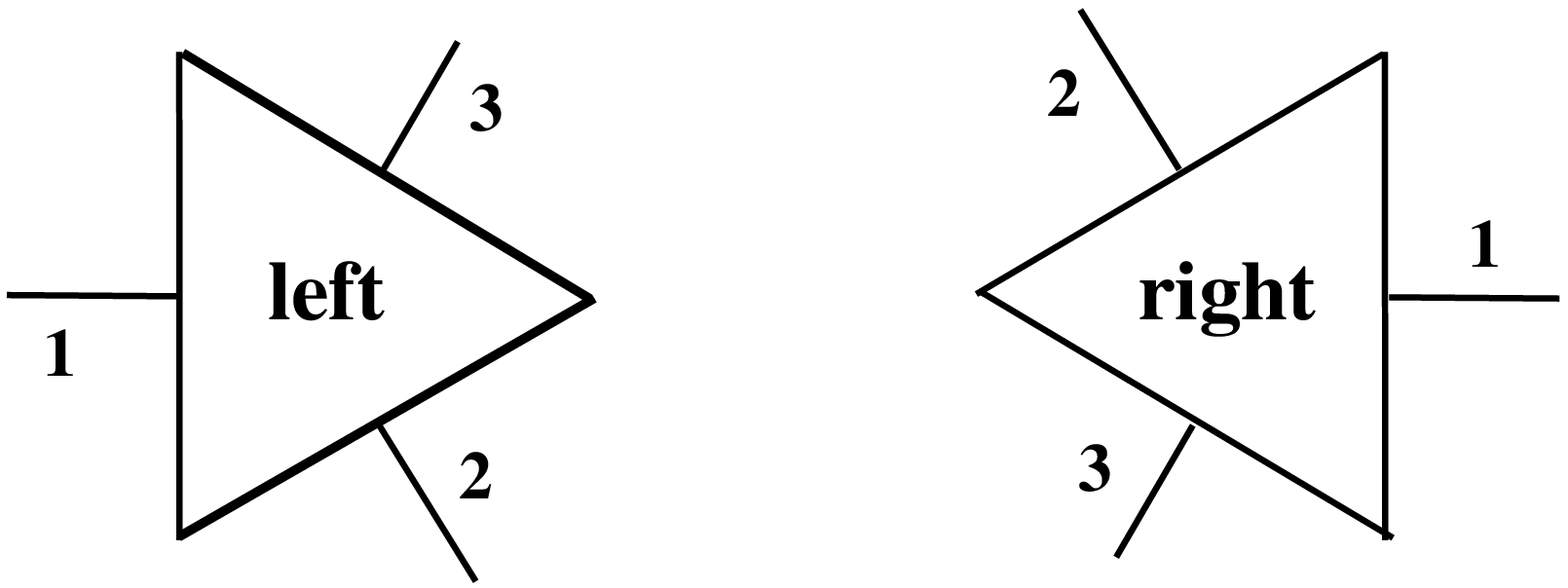} }   

\noindent  {\bf Figure 15}. \quad 
Two types of triangles for the D2T4 scheme with  equilaterals.   
The local numbers are explicited for each edge.  
\smallskip \smallskip 

  \smallskip   \noindent 
We observe that such a lattice contains only two types of equilateral triangles: 
the ``left'' and ``right''' triangles as displayed in Figure~15. 
We precise now the choices we have done to construct our scheme. 
The family $ \, {\cal P} \,$ of polynomials is simply composed by 
a restriction of (\ref{polynomes}) to the first four 
terms: 
$ \,\,   {\cal P} = \big\{  1, \, X ,\,  Y ,\, X^2 + Y^2    \big\} . \, $ 
Because we have two generic triangles, we have 
two families of neighboring directions 
$\,\,  \xi^{\rm left} = \big(-1, \, 0 \big) ,\, \big( {1\over2} ,\, -{{\sqrt{3}}\over{2}} \big)
 ,\, \big( {1\over2} ,\, {{\sqrt{3}}\over{2}} \big) ,$
$\,\, \xi^{\rm right} = \big(1, \, 0 \big) ,\, \big( -{1\over2} ,\, {{\sqrt{3}}\over{2}} \big)
 ,\, \big( -{1\over2} ,\, -{{\sqrt{3}}\over{2}} \big) .$  
We observe also that due to the simple numbering of local edges (see the figure~15), 
we have the simple relations 
$ \,\,  n_j(x) \equiv j \, $ and $ \,\,  \xi^{\rm left}_j +  \xi^{\rm right}_j = 0 .$ 
In this contribution, we consider only one 
 conserved variable  $ \, \, \rho =  m_0 \equiv \sum_{j=0}^3 f_j . \, $
The moments at    equilibrium are simply chosen with 
$ \,\,  m_1^{\rm eq} = m_2^{\rm eq} = 0 \,$ and 
$ \,  m_3^{\rm eq} = a_3 \, \rho  . \, $ 
%

\smallskip \noindent  {\bf Proposition 2. \quad Transition matrix for the D2T4 lattice
  Boltzmann scheme}

\noindent 
For the D2T4 lattice Boltzmann scheme defined previously, 
we have 
\moneq   \label{d2t4-PMm1}   
P^{\rm left}   \,=\, \big( M^{\rm left} \big)^{-1} \,, \qquad 
P^{\rm right}   \,=\, \big( M^{\rm right} \big)^{-1} \,.  
\monend   
In this particular case, the relations (\ref{d2t4-PMm1}) are exactly analogous to the ones
for  lattice Boltzmann schemes on Bravais lattices. 
In some sense, for the D2T4 scheme, the relations (\ref{d2t4-PMm1}) 
  remain (too !) simple ! 
%

  \smallskip  \noindent  {\bf Proof of Proposition 2}

 \noindent
Recall that due to (\ref{matrice-P}), we have 
  $ \,  P(x)_{i \,  \ell} = 
  \big(  \widetilde{M}(x_i) \big)^{-1} _{\,\,\,   n_i(x) \, \ell}  \, $ 
 with   $ \, n_i(x) \equiv i  \, $  due to our precise choice of numbering
(see the figure~15). Then we have the two matrix equalities  
$\, \,    P^{\rm left} _{\,\,\,i \,  \ell}  =
  \big(  \widetilde{M}^{\rm right} \big)^{-1} _{\,\,\, i  \, \ell}  \, \,$
and 
$ \,\,     P^{\rm right} _{\,\,\,i \,  \ell}  =
  \big(  \widetilde{M}^{\rm left} \big)^{-1} _{\,\,\, i  \, \ell} . $ 
We remark also    that 
  $ \, M^{\rm left}_{kj} = p_k(\xi^{\rm left}_j) \, $ and 
 $ \, \widetilde{M}^{\rm left}_{kj} = p_k(-\xi^{\rm left}_j) . \, $ 
Analogously   
 $ \, M^{\rm right}_{kj} = p_k(\xi^{\rm right}_j) \, $ and 
 $ \, \widetilde{M}^{\rm right}_{kj} = p_k(-\xi^{\rm right}_j) \, . $  
But     $ \,  \xi^{\rm left}_j +  \xi^{\rm right}_j = 0 \,, $ 
 then  
$ \,     P^{\rm left} = 
 \big(  \widetilde{M}^{\rm right} \big)^{-1} =   ( M^{\rm left} )^{-1} \, $  
and for the other family of triangles 
$ \,     P^{\rm right} = 
 \big(  \widetilde{M}^{\rm left} \big)^{-1} =   ( M^{\rm right} ) ^{-1} \, . $
The relation (\ref{d2t4-PMm1}) is established. \hfill $\square$


\smallskip   \monitem    Taylor expansion analysis for the D2T4 scheme  

\noindent  
The analysis  can   now be conducted without difficulty 
in the same framework than previously. We adopt the 
diffusive-scaling (\ref{diffusive-scaling}). 
After some developments with the help of formal calculus
(see {\it e.g.} \cite{DL09}) we derive the  
equivalent partial differential equation at the order 6:
\moneq   \label{d2t4-chaleur-o6} \left\{ \begin{array} [c]{rl}  
\displaystyle   {{\partial \rho}\over{\partial t}} \,-\, \mu \, \Delta \rho 
\, & = \,\,  \displaystyle  
{{a_3 \, \zeta}\over{24}} \, (12 \,\sigma_1^2 - 1) \, \Delta x \,  
\big( \partial_x^2 - 3 \,\partial_y^2 \big) \big( \partial_x \rho \big) 
+ \Theta_2 \, \Delta x^2 \,  \Delta^2  \rho   \\
\, &  \,\, + \quad    \Theta_3 \, \Delta x^3 \, 
\big( \partial_x^2 - 3 \,\partial_y^2  \big)   \Delta    \big( \partial_x \rho \big)   + 
\, \Delta x^4 \, A_6 \, \rho +   {\rm O}(\Delta x^6)  \, .  
\end{array} \right.     \monend
The notation $\, \sigma_k \,$ is identical to the one used at  
 H\'enon's relation \cite{He87}. The  diffusion coefficient 
$ \, \mu \, $  satisfies the relation
$\,\,   \mu =    \zeta \, a_3 \, \sigma_1  . \, $ 
%

  \smallskip \monitem  
``First order'', ``second order'', ``third order'' and ``quartic''  coefficients 

\smallskip \noindent
We have chosen $ \, \zeta = 1  . \,$ 
For first order simulations, 
we have taken    the following numerical values  
\moneq   \label{d2t4-order1}
a_3 =   0.216506350946109  ,\,\, 
s_1  =  1.2  ,\,\, 
s_3  =  0.750796078775233   \monend 
%
%
compatible  
with a diffusion coefficient 
$  \,\, \, \mu = 0.0721687836487032 = {{1}\over{4 \, \sqrt{12}}}  $. 
With the choice  $ \,\, \sigma_1 = {{1}\over{\sqrt{12}}} \,\, $ 
the scheme is at least second order accurate
(see the right hand side of (\ref{d2t4-chaleur-o6}))
and we take parameters 
to fit the previous choice of the diffusion coefficient: 
\moneq   \label{d2t4-order-sup2}
a_3 =   0.25   ,\,\,  
s_1  =  1.267949192431122     \monend 
%
%
With the particular value
\moneq   \label{d2t4-order2} 
s_3 = 0.422649730810374  \,,  \monend
we have $ \,   \Theta_2   \neq 0  \, $ 
and the D2T4 scheme is formally second order accurate. 
With 
\moneq   \label{d2t4-order3} 
s_3 = 0.758775495823486   \,,  \monend
we have $ \,   \Theta_2   = 0  , $ $  \,  \Theta_3   \neq 0 \, $
and the D2T4 scheme is formally third order accurate. 
With the choice of parameters 
\moneq   \label{d2t4-order4}  
s_3 = 0.732050807568877   = \sqrt{3} - 1  \,,  \monend
{\it id est} $ \,  \sigma_3 = {{\sqrt{3}}\over{2}}  , \,$ 
we have $ \,   \Theta_2   =   \Theta_3   = 0 . \, $
With these conditions, the D2T4 scheme is theoretically fourth order accurate.

  \newpage 
\bigskip \bigskip   \noindent {\bf \large 7) \quad  Diffusion simulations with the D2T4 scheme}  


\noindent    
We have done essentially the same simulations  as performed  with the D2T7 lattice
Boltzmann scheme (see Section 4). 
 
 \smallskip  \monitem     One point periodic analysis     
  
\noindent 
The results are presented in Figure~16. The theoretical orders with the four 
choices
of parameters proposed previously are exactly the one proposed by the Taylor expansion
analysis. 
A  
defect  of isotropy for the numerical diffusivity  is clearly visible for
parameters that lead 
to a first order and third order schemes  with this 
 D2T4 simulator.

\smallskip                        
\centerline { \includegraphics[width=.47 \textwidth] {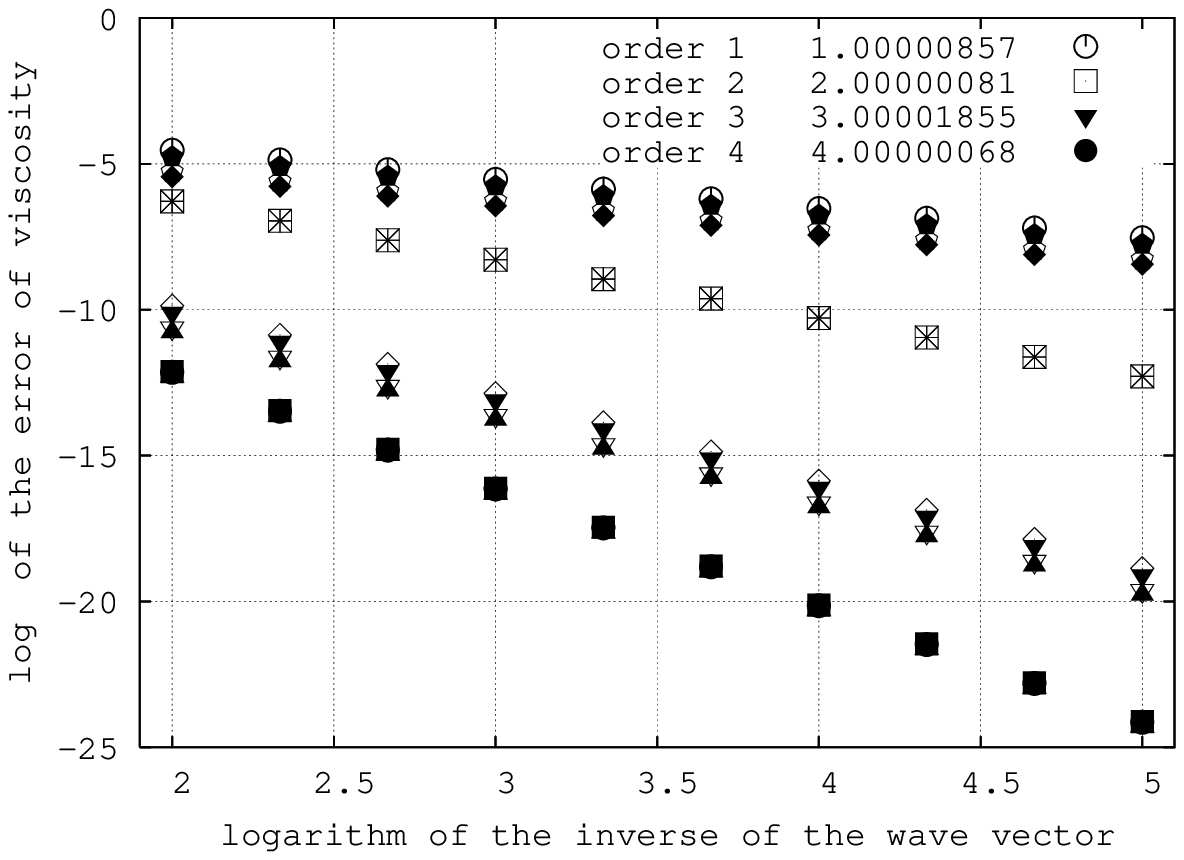}} 

\noindent  {\bf Figure 16}. \quad 
One point periodic analysis. 
Error  $ \epsilon \, \equiv \,   \mid \mu -  \mu_{\rm num} \mid \, $  
 between  numerical and theoretical diffusivities. Four sets of parameters 
defined at relations (\ref{d2t4-order1}), (\ref{d2t4-order-sup2}), 
  (\ref{d2t4-order2}),  (\ref{d2t4-order3}) and   (\ref{d2t4-order4})
lead to schemes of several orders. The measured orders with a linear regression 
are displayed in the right column.  
\smallskip \smallskip 

\smallskip    \monitem       Periodic pipe and rectangle 

\noindent    
We have tested the fourth order version  (\ref{d2t4-order-sup2}) (\ref{d2t4-order4})
of the D2T4 lattice Boltzmann scheme on two simple periodic geometries presented in
Section~4. The numerical results (Figure~17) show that the scheme
is convergent, but simply at second order accuracy.

\smallskip                     
 \centerline { \includegraphics [width=.47 \textwidth] {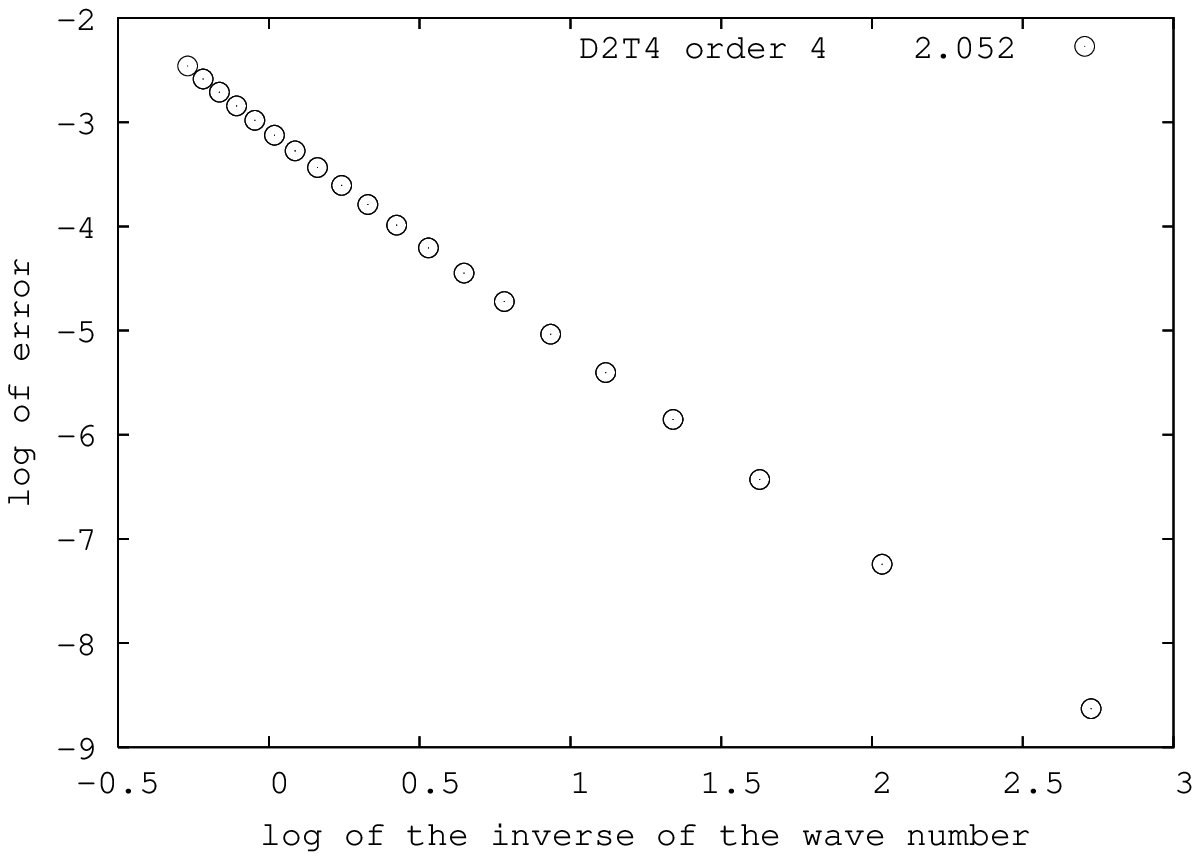} 
\quad 
 \includegraphics[width=.47 \textwidth] {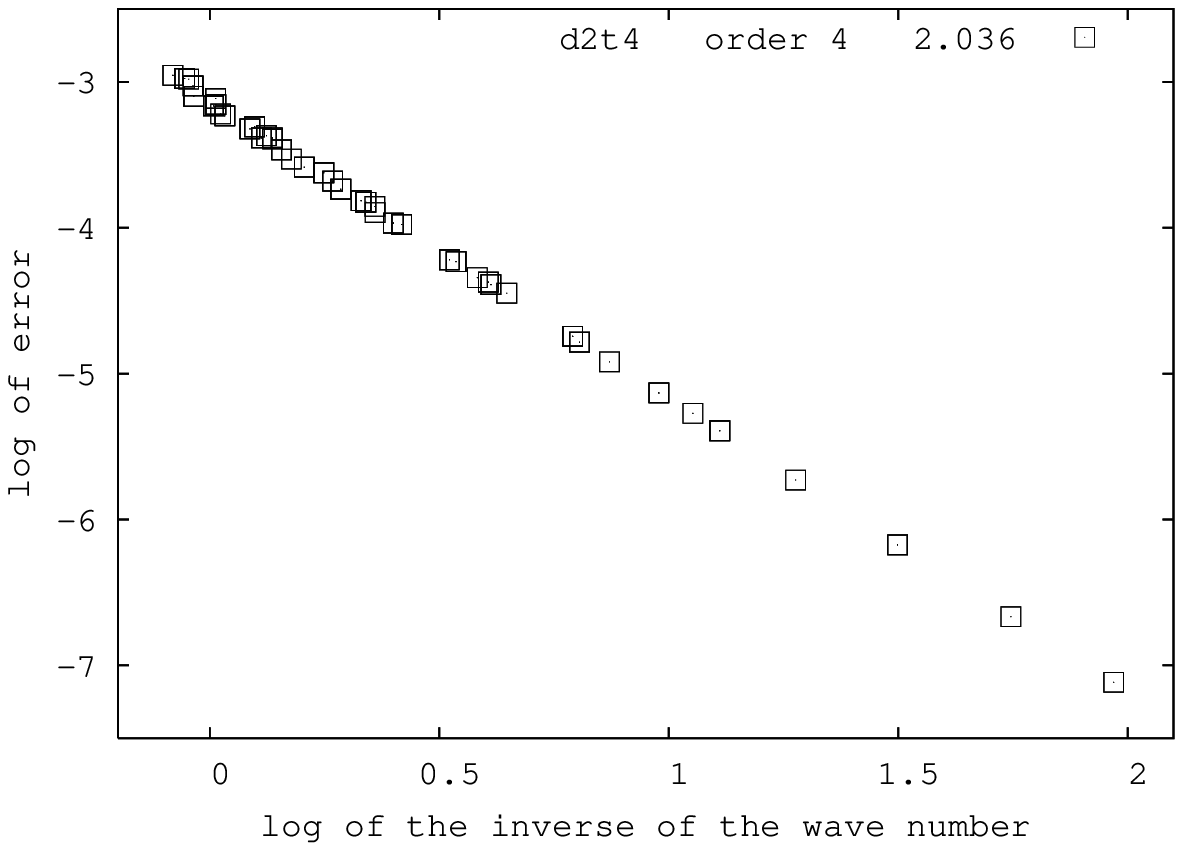}  }

\noindent  {\bf Figure 17}. \quad 
D2T4  lattice Boltzmann scheme for the heat equation. 
Periodic modes for a pipe with $ nx = 96 $  and $ ny = 4 $ mesh points (left)
and   a rectangle of  $ nx = 36 $  by $ ny = 52 $     points (right). 
Error $ \epsilon \, \equiv \,   \mid \mu -  \mu_{\rm num}  \mid \, $  
between  exact and numerical diffusivities. The predicted coefficients for the 
order 4 define a second order scheme in this particular case.    
\smallskip  \smallskip  

\smallskip   \monitem 
Harmonic polynomials on a triangle     

\noindent  
The numerical computation of Laplace equation with non-homogeneous boundary conditions 
has been also performed by integrating  the heat equation and taking
the limit for time large enough.
Our simulation   (Figure~18) shows that the D2T4 scheme is convergent
with second order accuracy. Nevertheless the second order and fourth 
order versions of the scheme give essentially the same results.

\smallskip                      
 \centerline {  \includegraphics [width=.47 \textwidth] {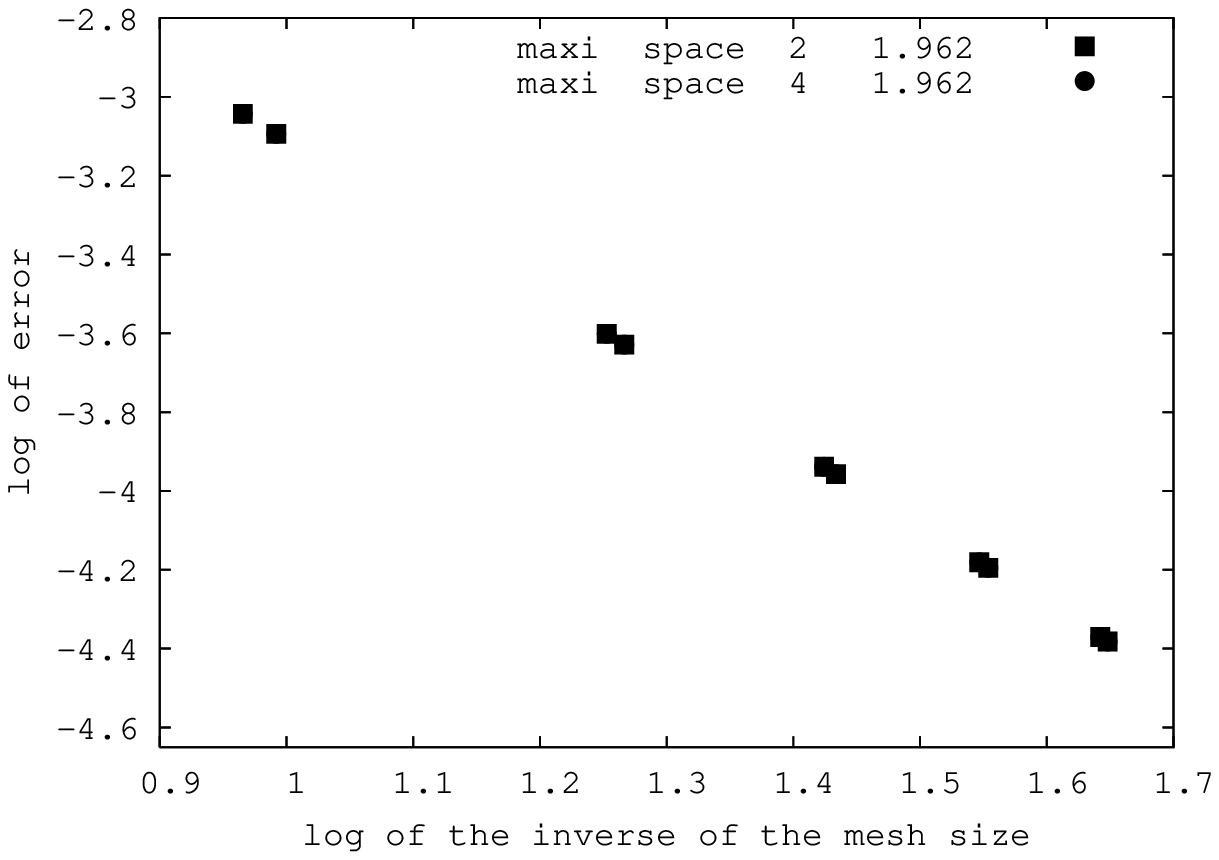} 
\quad  
\includegraphics [width=.47 \textwidth]  {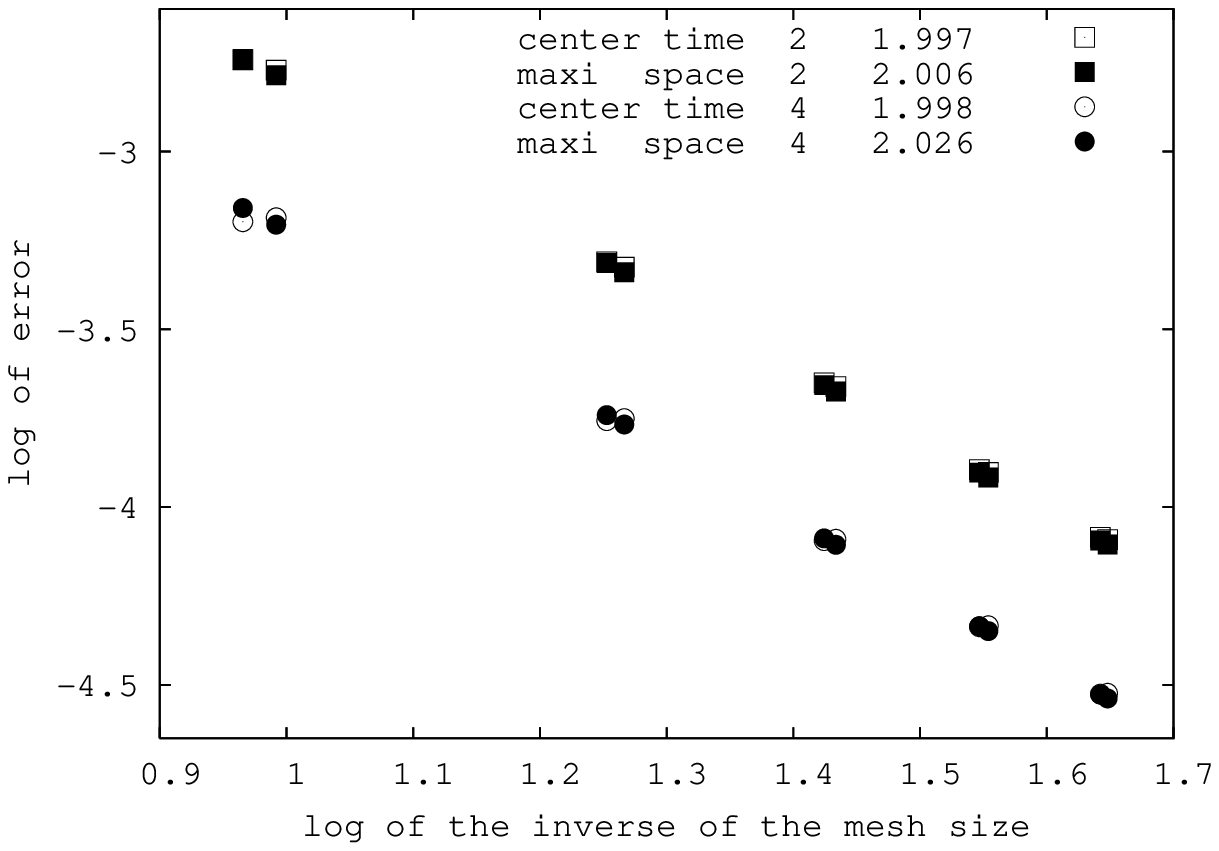}  } 

 \noindent  {\bf Figure 18}. \quad 
Two-dimensional computation of the harmonic function  $ \, p_H(x,\, y)  = x^2 - y^2 \, $  
on a triangle with the D2T4 lattice Boltzmann scheme (left).  
Convergence of the   $ {\rm L}^\infty $ error  for several meshes. 
The D2T4 lattice Boltzmann scheme remains  of order 2 even if quartic parameters are 
used in the simulation.   
  Dissipation of the first Dirichlet mode for an equilateral triangle (right).  
   Time and $ {\rm L}^\infty $ space errors 
for several meshes and several ``orders''. 
The  numerical accuracy is equal to 2 for all the parameters. With ``quartic'' parameters
the absolute level of the error is substantially reduced.    
\smallskip 
 
\smallskip   \monitem      Dissipation of a triangular Dirichlet mode

\noindent 
The dissipation of the first mode  described at the figure~7 
for the D2T7 scheme has been constructed without difficulty.
Now the two main versions of the scheme (second and fourth orders)
converge with second order accuracy as shown in Figure~18. 
We observe that even  if no extra order of convergence has been obtained, 
the results with quartic parameters give a better precision.

\smallskip    \monitem      Dirichlet modes for a triangle  

\noindent 
The simulations  done with the D2T7 lattice Boltzmann scheme have been 
compared with a D2T4 simulator. 
The results (Figure~19) explicit this comparison. 
A first result is that the level of error for D2T4 is comparable with D2T7 results at
order two.
If we look precisely to the error fields, 
distinguished contribution is due to the boundary conditions. 
    
 \smallskip
                          
\noindent 
 \includegraphics[width=.37 \textwidth, angle=90] {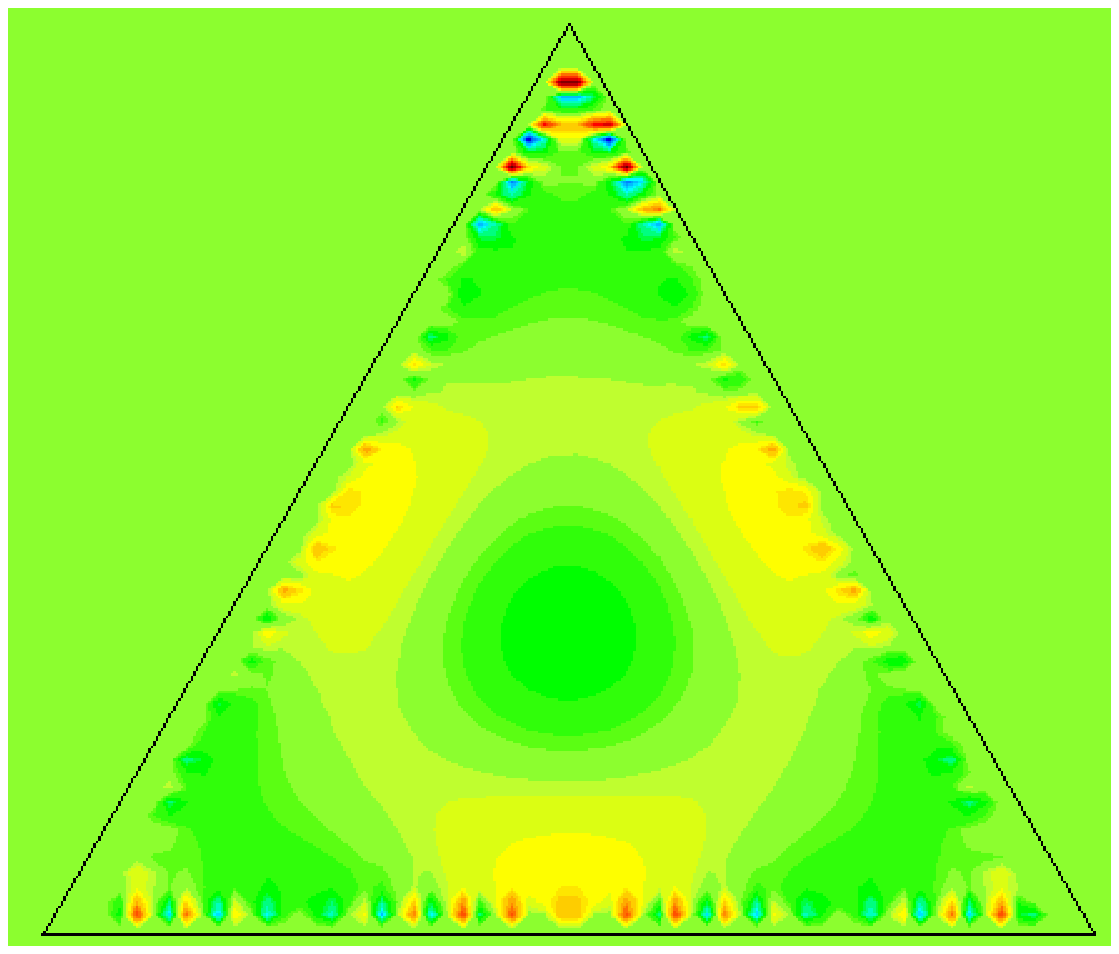} $\,\,$
 \includegraphics[width=.37 \textwidth, angle=90] {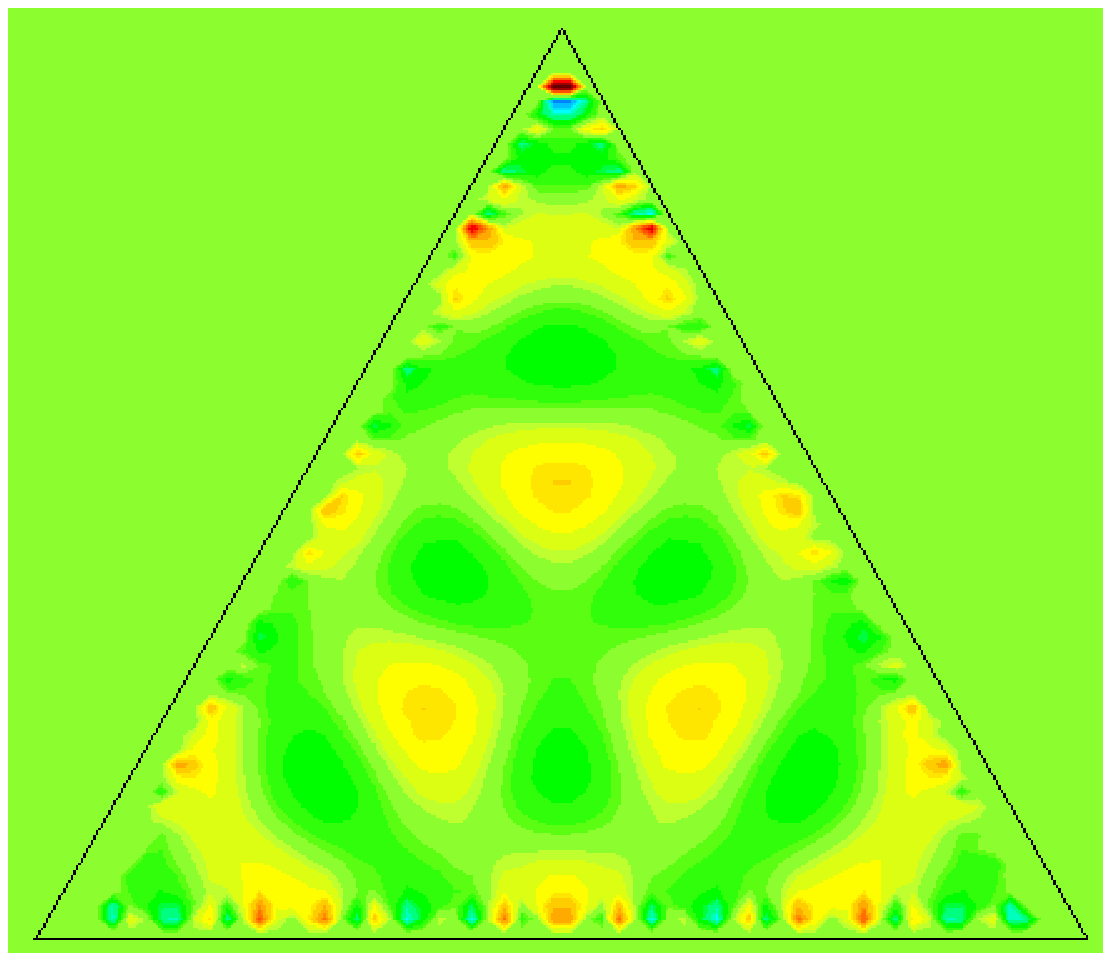} $\,\,$
 \includegraphics[width=.37 \textwidth, angle=90] {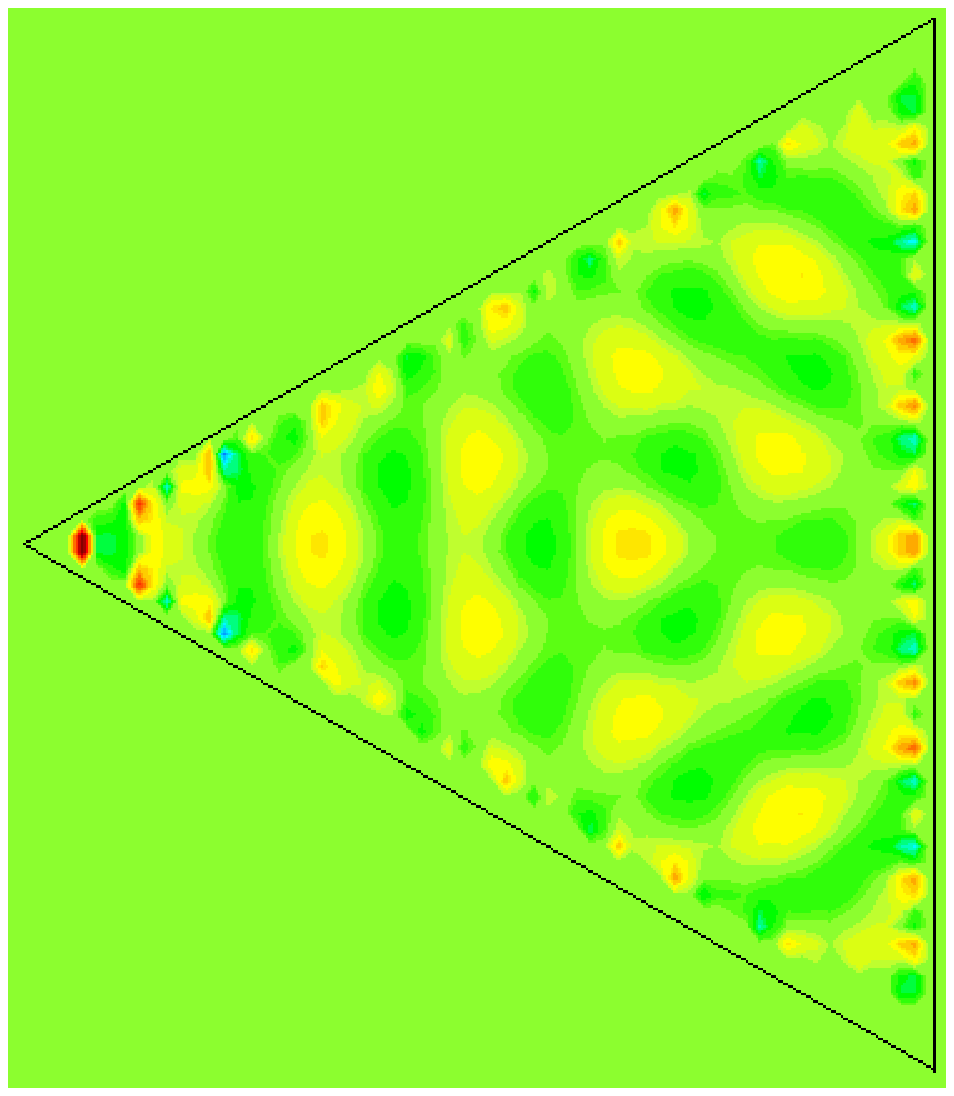}  

\noindent  {\bf Figure 19}. \quad 
Isovalues of the errors for D2T4 computation of     Dirichlet modes  of an equilateral triangle. 
 The exact reference eigenvalues are equal respectively to 12, 48 and 108
and the   computed ones to 11.97493,    47.59816  and  105.95870. 
This figure can be compared to second order accurate D2T7 results at 
Figures~9, 10 and 11 respectively.  
\smallskip \smallskip 


\bigskip \bigskip   \noindent {\bf \large 8) \quad  Conclusion}  
 

\noindent   
We have proposed an extension of the lattice Boltzmann method for triangular
meshes. Our first step concerns a single conservation law and we 
made numerical simulations  for the heat equation.
For an extension of the discrete particle method, we have considered non Bravais
lattices such that the degree of  each vertex is constant.  
Our formulation does not need any  finite volume  or Delaunay-Voronoi triangulation
hypothesis as in the previous contributions. 
We have used the Taylor expansion analysis with a diffusive scaling to explicit some
parameters of the d'Humi\`eres scheme. With this method, it is possible to get formally a
better accuracy. Our simulations show that this extra accuracy can be 
 obtained with very fundamental one point periodic hypothesis. In  more  
realistic cases, this extra-accuracy is in general not   observed.

We think that triangular meshes explicit the limit of validity of the Taylor expansion
analysis. 
In fact when we write the lattice Boltzmann scheme with the relation  (\ref{LB-scheme})
or (\ref{LB-particle}) and when we perform the Taylor expansion, 
we suppose that there exists a very regular function $ \, f(x, \, t) \,$ 
of space and time that support the definition of the scheme.
In particular, this function is supposed to be independent of the lattice~!
This last Ansatz is   in defect  for triangular meshes on nonsymmetric lattices 
as D2T4. Note that this kind of remark recover other critics \cite{Ch90, GS86}
relative to this kind of symptotic analysis  \cite{LP74, WH74}.

Two directions of research are natural in the continuation of the present contribution. 
First we can try to develop a true mathematical analysis of the lattice Boltzmann 
scheme, following {\it e.g.} previous work of Junk and Yong \cite {JY09}    
with appropriate mathematical tools,  as done typically by Ciarlet  and Raviart 
for finite elements \cite{CR72} or Gallou\"et and coworkers  for finite volumes 
\cite{CGH93}. 
Second we can extend triangular lattice Boltzmann schemes to systems 
with other conservation laws for acoustics and fluid flow applications, 
revisiting the breakthrough of Frisch,  Hasslacher and  Pomeau \cite {fhp}. 
Preliminary results  have been obtained for D2T10, which are not
described in this article  due to unnecessary complications.

\bigskip   \bigskip  \noindent {\bf \large Acknowledgments}   


 \noindent    
The authors  thank 
 the ``LaBS project'' (Lattice Boltzmann Solver,  www.labs-project.org), 
funded by the French ``FUI8 research program'', for supporting this contribution. 
Last but not least, 
the authors thank the referees for  very constructive remarks. 
Some of them  have been   
incorporated into  the present edition of the article.


\medskip

\end{document}